\documentclass{article}
\usepackage{latexsym}
\usepackage{amsfonts}
\usepackage{amsmath}
\usepackage{amsfonts}
\usepackage{epsfig}
\usepackage{stmaryrd}
\usepackage{amssymb}
\usepackage{pxfonts}
\usepackage{wasysym}
\usepackage[font={small}]{caption}
\usepackage{algorithm}
\usepackage{algorithmic}
\newtheorem{theorem}{Theorem}
\usepackage{algorithm}
\usepackage{algorithmic}
\usepackage{bm}
\usepackage{bbm}

\newtheorem{definition}[theorem]{Definition}
\newtheorem{lemma}[theorem]{Lemma}

\title{A multigrid-like algorithm for probabilistic domain decomposition}
\date{\vspace{-5ex}}
 \author{ Francisco Bernal \footnotemark[1] \ \footnotemark[3]
\and Juan A. Acebr\'{o}n \footnotemark[2]\ \footnotemark[1]}

\begin{document}
\maketitle

\renewcommand{\thefootnote}{\fnsymbol{footnote}}

\footnotetext[1]{INESC-ID$\backslash$IST, TU Lisbon. Rua Alves Redol 9, 1000-029 Lisbon, Portugal.}
\footnotetext[2]{ISCTE - Instituto Universit\'ario de Lisboa
Departamento de Ci\^{e}ncias e Tecnologias de Informa\c{c}\~{a}o. Av. das For\c{c}as Armadas 1649-026 Lisbon, Portugal. ({\tt juan.acebron@ist.utl.pt})}
\footnotetext[3]{Center for Mathematics and its Applications,
Department of Mathematics, Instituto Superior T\'ecnico. 
Av. Rovisco Pais 1049-001 Lisbon, Portugal. ({\tt francisco.bernal@ist.utl.pt})}
\renewcommand{\thefootnote}{\arabic{footnote}}

\begin{abstract}
We present an iterative scheme, reminiscent of the Multigrid method, to solve large boundary value problems with Probabilistic Domain Decomposition (PDD). In it, increasingly accurate approximations to the solution are used as control variates in order to reduce the Monte Carlo error of the following iterates--resulting in an overall acceleration of PDD for a given error tolerance. The key ingredient of the proposed algorithm is the ability to approximately predict the speedup with little computational overhead and in parallel. Besides, the theoretical framework allows to explore other aspects of PDD, such as stability. One numerical example is worked out, yielding an improvement of between one and two orders of magnitude over the previous version of PDD. 
\end{abstract}


\section{Introduction} \label{S:Introduction}
 
{\bf Probabilistic and deterministic domain decomposition.} In the solution of large boundary value problems (BVPs) arising in realistic applications
, the discretization of the BVP on a domain $\Omega$ 
leads to algebraic systems of equations that can only be solved on a parallel computer with $P\gg 1$ processors by means of domain decomposition. The idea is to divide $\Omega$ into a set of $m\geq P$ overlapping subdomains, $\Omega= \cup_{i=1}^m \Omega_i$, and have processor $j$ solve the restriction of the partial differential equation (PDE) to the subdomain $\Omega_j$. Not only does parallellization need parallel computers but parallel algorithms as well. State-of-the-art methods---which we will refer to as 'deterministic'--are iterative and require updating on the interfaces \cite{Smith96}, a step which unavoidably involves interprocessor communication and thus is intrinsically sequential. Regardless of the sophistication of the deterministic method, this will eventually set an upper limit to the scalability of the algorithm according to Amdahl's law. Whether or not this is a practical concern depends on the size of the problem, or, equivalently: the size of the problems which can be tackled is ultimately determined by the scalability limit of the domain decomposition algorithm.

An alternative to deterministic methods which is specifically designed to circumvent the scalability issue is the probabilistic domain decomposition (PDD) method, which has been successfully applied to elliptic \cite{Acebron2005} and parabolic BVPs \cite{Acebron2009}\cite{Acebron2013}. 
PDD consists of two stages. 
In the first stage, the solution is calculated only on a set of interfacial nodes along the fictitious interfaces, by solving the probabilistic representation of the BVP with the Monte Carlo method. More precisely, the pointwise solution of the BVP is $u(t,{\bf x})=E[\phi({\bf X}_s)|{\bf X}_t={\bf x}]$, i.e. the expected value of a functional $\phi$ of a given stochastic process ${\bf X}_s$ conditioned to $(t,{\bf x})$. It is then possible to reconstruct (approximately) the solution on the interfaces, so that the PDE restricted to each of the subdomains is now well posed, and can be independently solved--the second stage of PDD. 
Note that both stages in PDD are embarrasingly parallel by construction. Moreover, PDD is naturally fault-tolerant.
\newline

\begin{figure}[h]
\centerline{\includegraphics[width=1.\columnwidth]{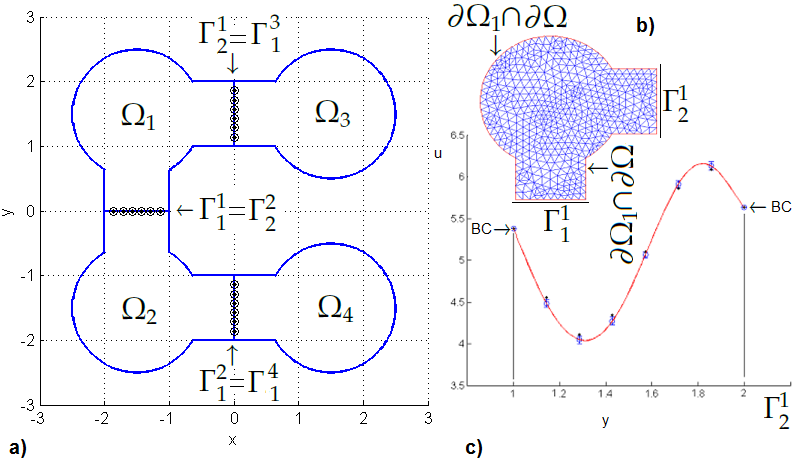}}
\caption{{\bf a)} PDD partition into $4$ subdomains and $3$ interfaces with $6$ nodes each ($n=18$). 
{\bf b)} FEM mesh on subdomain $\Omega_1$. {\bf c)} The nodal values interpolated with RBFs along the interface $\Gamma_2^1=\Gamma_1^3$, making up a Dirichlet BC for $\Omega_1$ and $\Omega_3$.}
\label{I:Figura1}
\end{figure}

{\bf Nomenclature of PDD.} In this paper, we shall exclusively deal with elliptic BVPs. Let us introduce some terminology (see Figure \ref{I:Figura1}). The domain $\Omega\subset{\mathbb R}^D,\,D\geq 2$ (not necessarily simply connected) on which the BVP is being solved is partitioned into $m$ {\em nonoverlapping} subdomains $\Omega_1,\ldots,\Omega_m$. 
The boundary $\partial\Omega_k$ of a subdomain $\Omega_k$ contains several ($m_k\geq 1$) {\em artificial interfaces}--each of which is shared between $\Omega_k$ and another adjacent subdomain--which are labeled $\Gamma_1^k,\ldots,\Gamma_{m_k}^k$ (note that this labeling is not unique). 
A subdomain boundary $\partial\Omega_k$ may or not contain some portion of the actual boundary. 
In sum, $\partial\Omega_k= \big(\partial\Omega_k\cap\partial\Omega\big) \,\cup\, \big(\cup_{j=1}^{m_k}\Gamma_j^k\big)$. Artificial interfaces are discretized into {\em interfacial nodes} (or simply, nodes) uniquely labeled $\{{\bf x}_1,\ldots,{\bf x}_n\}$. 
\begin{definition}Assume $m$ functions $f_i({\bf x})$ defined on $m$ domains $D_i$ such that $f_i({\bf y})=f_j({\bf y})$ if ${\bf y}\in D_i\cap D_j$. Then, their direct sum is defined as
\begin{equation}
\oplus_{k=1}^m f_k({\bf x})=\left\{
\begin{array}{ll}
f_i(\bf x) & \textrm{ if }{\bf x}\in D_i\, (i=1,\ldots,m),\\
0 & \textrm{ if }{\bf x}\notin \cup_{1\leq k \leq m}D_k.\\
\end{array}
\right.
\end{equation}
\end{definition}
The solution of the BVP on the nodes is calculated by resorting to the probabilistic formulation of the BVP with a Monte Carlo method, yielding the {\em nodal values} (or {\em nodal solutions}) $\{u_1,\ldots,u_n\}$. 
Consider a subdomain $\Omega_k$. A Dirichlet BC can then be provided on every $\Gamma_j^k$ by interpolation of the nodal values $\{u_i\,|\,{\bf x}_i\in\Gamma_j^k\}$. Along with the actual BCs which apply on $\partial\Omega_k\cap\partial\Omega$, the BVP on $\Omega_k$ now is well posed and can be solved right away, yielding $v_k({\bf x})$. Once the subdomain solutions $v_1({\bf x}),\ldots,v_m({\bf x})$ are available, they are put together to form a global PDD solution: ${\tilde u}({\bf x})=\oplus_{k=1}^m v_k({\bf x})$. (We reserve $u({\bf x})$ for the exact solution of the BVP and denote the global PDD approximations with ${\tilde u}({\bf x})$.)
Since adjacent subdomains $\Omega_i$ and $\Omega_j$ share a Dirichlet BC on their common interface, the PDD solution 
is continuous in $\Omega$--although not necessarily differentiable.\newline 


{\bf The cost of computing the nodal solutions.} The convenience of PDD depends on whether a suitable stochastic representation for the BVP under consideration is available, and on the cost involved in numerically solving it. Due to the poor accuracy of the Monte Carlo method (compared with deterministic ones), the bulk of the cost of PDD falls on the calculation of the nodal solutions to within a required accuracy. More precisely, given a {\em nodal error tolerance} $a$ and a {\em confidence interval} $P_q$, the cost of solving the BVP on an interfacial node scales as ${\cal O}(a^{-2-1/\delta})$, where $\delta$ is the weak convergence rate of the numerical integration scheme (also called the {\em bias}). This poor rate of convergence is due to the slow convergence of both the statistical error and the bias \cite{Kloeder&Platten99}, which have to be tackled simultaneously. For BVPs with Dirichlet BCs there quite a few linear integrators (i.e. with $\delta=1$) \cite{OurSurvey,Gobet&Menozzi_2010,Milstein&Tretyakov_book}. 

Regarding the statistical error, replacing the mean by a Multilevel estimator of the expected value of Feynman-Kac functionals has recently been shown to dramatically reduce the cost to ${\cal O}(|\log a|^3a^{-2})$ \cite{Giles_y_yo}. When the bias law is sharp, extrapolation \'a la Talay-Tubaro or regression methods in the spirit of \cite{LTMCR} can further improve the accuracy at virtually no extra cost.\newline 

{\bf Using rougher numerical solutions to reduce the variance.} By construction, PDD offers an additional device to accelerate the Monte Carlo simulation of the nodal values, namely the possibility of calculating and exploiting rougher estimates of the global solution of the BVP. 
Assuming that a numerical solution ${\tilde u}_0({\bf x})$ with a {\em target nodal error tolerance} $a_0$ is required, it may be worth to calculate before a rougher approximation ${\tilde u}_1({\bf x})$, with an $a_1>a_0$ tolerance; and then use it to draw the stochastic pathwise nodal control variate 
\begin{equation}
\label{F:Intro_xi}
\xi= -\int_0^{\tau}e^{\int_0^t c({\bf X}_s)ds}\sigma^T\nabla{\tilde u}_1\,d{\bf W}_t,
\end{equation}    
alongside the Monte Carlo realizations of the Feynman-Kac functional. (In (\ref{F:Intro_xi}), the integrals are Ito's, $c$ is the BVP potential, and $\sigma$ and $\tau$ are the diffusion matrix and first-exit time from $\Omega$ of the stochastic process (\ref{F:SDE}) driven by a Wiener process ${\bf W}_t$.) 
This allows one to construct afterwards an estimator of the Feynman-Kac functional involving the control variate (\ref{F:Intro_xi}), which has the same expected value but a much smaller variance. This notion is what we call IterPDD. In order to fix ideas, let us 
introduce the following notation:
\begin{equation}
\label{F:PlainPDD}
{\tilde u}(a)= \textrm{PlainPDD$(a)$, or simply PlainPDD$(a)$}
\end{equation}
means that ${\tilde u}({\bf x})$ is a PDD approximation obtained with tolerance $a$ and no variance reduction; while
\begin{equation}
\label{F:IterPDD_2}
[{\tilde u}_0(a),\xi_0({\tilde u}_1)]= \textrm{IterPDD$(a_0,a_1)$, or simply IterPDD$(a_0,a_1)$} 
\end{equation}
means that first ${\tilde u}_1(a_1)= \textrm{PlainPDD}(a_1)$ is calculated without variance reduction, then differentiated in order to construct $\xi$ according to (\ref{F:Intro_xi}), which in turn is used as control variate in order to reduce the variance in calculating ${\tilde u}_0$ with a target tolerance $a_0$, which is the ultimate goal. Because the nodal values of ${\tilde u}_0$ can now be calculated with much less variance, statistical errors are smaller, and the time (or cost) it takes the computer to hit the tolerance $a_0$ is also less. In fact,
\begin{eqnarray}
\label{F:Cost_of_IterPDD}
\textrm{cost of IterPDD$(a_0,a_1) \approx$
cost of PlainPDD$(a_1)$} + \nonumber\\
\textrm{$\big(1 - \rho^2[\phi_0,\xi_0({\tilde u}_1)]\big)\times$cost PlainPDD$(a_0)$},
\end{eqnarray}

where $\rho[\phi_0,\xi_0({\tilde u}_1)]$ is Pearson's correlation. 
As (\ref{F:Cost_of_IterPDD}) indicates, there is a tradeoff between the effort invested in calculating ${\tilde u}_1(a_1)$, and the reduction of variance yielded by $\xi_0({\tilde u}_1)$, which depends on the quality of ${\tilde u}_1(a_1)$. The most straightforward procedure is to simply guess some $a_1>a_0$. While numerical tests indicate that the IterPDD strategy can be quite successful, a poorly chosen $a_1$ may well result in an overall cost of IterPDD$(a_0,a_1)$ {\em larger} than that of PlainPDD$(a_0)$. 
Therefore, at the heart of IterPDD lies an optimization problem for $a_1$. In order to make educated guesses of $a_1$ given $a_0$, two questions must be tackled: i) how does the cost of a ${\tilde u}(a)=$PlainPDD$(a)$ simulation depend on $a$; and ii) how does $\rho[\phi_0,\xi_0\big({\tilde u}(a)\big)]$ depend on $a$. The former requires that the SDE integrator have a predictable and sharp order of weak convergence.  Deriving a {\em sensitivity formula} for ii) is one of the main points of this paper.\newline

{\bf A multigrid-like PDD algorithm.} 
With such a sensitivity formula in place which can predict an optimal (or more realistically, good enough) initial tolerance $a_1$ for $[{\tilde u}_0(a_0),\xi_1({\tilde u}_1)]=$IterPDD$(a_0,a_1)$, it is natural to try and compute ${\tilde u}_1(a_1)$ faster by finding $a_2>a_1$ which minimizes the cost of $[{\tilde u}_1(a_2),\xi_2({\tilde u}_2)]=$IterPDD$(a_1,a_2)$. Much like in the Multigrid method, a number $J$ and an {\em optimal sequence} of nested IterPDD simulations can be envisioned with tolerances $a_J>\ldots>a_2>a_1>a_0$ which fully exploits the potential of control variates for a given BVP. However, in IterPDD, given a target tolerance $a_0$, the number $J$ and the sequence $a_J>...>a_0$ must be determined {\em before} actually running one single PDD simulation. To compound matters, all of $J$, $a_J,\ldots,a_1$, and in general any result provided by any such {\em scheduling algorithm} will be affected by the randomness introduced by Monte Carlo, making its performance meaningful only in terms of its expected value and variance. \newline 


{\bf Outline of the paper.} We start by revisiting the probabilistic formulation of elliptic BVPs with Dirichlet BCs in Section \ref{S:Representation}, 
and identifying the
pathwise control variates. 
Section \ref{S:Theory} formally poses the problem--namely, acceleration of PDD with the IterPDD scheme--and presents our own theoretical results concerning the aforesaid sensitivity $\rho[\phi_0,\xi_0\big({\tilde u}(a)\big)]$. 
Section \ref{S:GlobalError} links the nodal target error tolerance, $a_0$, to the global PDD error tolerance $\epsilonup$
and discusses the stability of PDD. In Section \ref{S:Approximations}, the formal restrictions in Section \ref{S:Theory} are relaxed leading to practical, but partially heuristic, scheduling algorithm/sensitivity formula (Algorithm \ref{A:Scheduling}) and final iterative, multigrid-like PDD loop (Algorithm \ref{A:IterPDD}). They are numerically tested in Section \ref{S:Experiment}, and Section \ref{S:Conclusions} concludes the paper.\newline

\section{Pathwise control variates for elliptic BVPs}
\label{S:Representation}

\subsection{Probabilistic representation}\label{SS:Dynkin}


Consider the linear elliptic BVPs with Dirichlet BCs in $\Omega\subset{\mathbb R}^D$:


\begin{equation} 
\label{F:EllipticBVP} 
\left\{
\begin{array}{ll}
L({\bf x})u + c({\bf x})u = f({\bf x}), 
& \textrm{ if } {\bf x}\in\Omega,\\
u({\bf x})= g({\bf x}), & \textrm{ if } {\bf x}\in\partial\Omega,
\end{array}
\right.
\end{equation} 

where the {\em differential generator} $L$ is
\begin{equation}
L({\bf x}):=\frac{1}{2}\sum_{i=1}^D\sum_{j=1}^D a_{ij}({\bf x})\frac{\partial^2}{\partial x_i\partial x_j} + {\bf b}({\bf x})\cdot\nabla,
\end{equation} 

and the functions $a_{11},\ldots,a_{DD},b_1,\ldots,b_D,c,f,g:{\mathbb R}\mapsto{\mathbb R}$ are regular enough that the solution to (\ref{F:EllipticBVP}) exists and is unique \cite{Milstein&Tretyakov_book}. In particular, if $c\leq 0$ and the matrix $A({\bf x})=[a_{ij}({\bf x})]$ is such that 
\begin{equation}
\Lambda\in{\mathbb R}:=\sup\limits_{{\bf x}\in\Omega\backslash\partial\Omega} \textrm{spectrum of $A({\bf x})$}\geq
\inf\limits_{{\bf x}\in\Omega\backslash\partial\Omega} \textrm{spectrum of $A({\bf x})$}=:\lambda>0,
\end{equation}
the pointwise solution to (\ref{F:EllipticBVP}) admits the following probabilistic representation 

\begin{equation}
\label{F:Dynkin}
u({\bf x}_0)= E[\phi]:= E\Big[\, g({\bf X}_{\tau})e^{\int_0^{\tau}c\big({\bf X}_s\big)ds} \, - \, \int_0^{\tau}f\big({\bf X}_t\big)e^{\int_0^t c\big({\bf X}_s\big)ds}dt\,\Big],
\end{equation}

where $E[\cdot]$ stands for the expected value, the functional $\phi$ is called {\em score}, and ${\bf X}_{\tau}$ is the value at $t=\tau>0$ of the stochastic process ${\bf X}_t:[0,\tau]\rightarrow\mathbb{R}^D$, driven by the stochastic differential equation (SDE)

\begin{equation}
\label{F:SDE}
d{\bf X}_t= {\bf b}({\bf X}_t)dt + \sigma({\bf X}_t)d{\bf W} _t,\qquad{\bf X}_0={\bf x}_0,
\end{equation}

where the diffusion matrix $\sigma$ is obtained from $\sigma\sigma^T=A$ (by Choleski's decomposition), and ${\bf W}_t$ is a standard $D$-dimensional Wiener process. 
The {\em first exit time} $\tau$ is defined as $\tau= \inf_{t} {\bf X}_t\in\partial\Omega$, i.e. the time when a solution of the SDE (\ref{F:SDE}) first touches $\partial\Omega$ at the {\em first exit point} ${\bf X}_{\tau}$. The process ${\bf X}_t,\,0\leq t\leq \tau$ can be thought of as a non-differentiable trajectory inside $\Omega$. Equation (\ref{F:Dynkin}) is Dynkin's formula, a particular case of the more general Feynman-Kac formula for parabolic BVPs--see \cite{Freidlin_Book} for more details.





\subsection{Errors arising in a Monte Carlo simulation}\label{SS:Errors}

In practice, solving (\ref{F:Dynkin}) involves two levels of discretization. First, the SDE (\ref{F:SDE}) has to be integrated numerically according to a numerical scheme $\Xi$ (which we will call {\em integrator}) with a timestep $h>0$, yielding a discretized score $\phi_h$. This results in $E[\phi_h]$ being a biased estimator of $E[\phi]$, with {\em signed} bias

\begin{equation}
\label{F:Signed_bias}
B_{\Xi}\big(E[\phi_h]\big):=E[\phi_h]-E[\phi].
\end{equation}

As suggested by the notation, the bias depends on both the timestep $h$ and the specific integrator $\Xi$. It takes asymptotically the form of a power law \cite{Kloeder&Platten99}:
\begin{equation}
B_{\Xi}\big(E[\phi_h]\big)\shortrightarrow \beta h^{\delta} \textrm{ as }h\shortrightarrow 0^+,
\end{equation}

where $\beta$ is a signed constant independent of $h$. Second, the expected value in (\ref{F:Dynkin}) is replaced by an estimator over $N$ independent realizations $\phi_h^{(1)},\ldots,\phi_h^{(N)}$ of the SDE (\ref{F:SDE})--which is the essence of the Monte Carlo (MC) method. Typically, that estimator is the mean, which, according to the central limit theorem, introduces a {\em statistical error} ($V[\cdot]$ is the variance):

\begin{equation}
\label{F:Statistical_error}
\big| E[\phi_h] - \frac{1}{N}\sum_{j=1}^N \phi_h^{(j)}\big|\leq q\sqrt{\frac{V[\phi_h]}{N}}
\end{equation}

with probability $P_q=68.3\%, 95.5\%$, and $99.7\%,$ for $q=1,2,3$ \cite{Kloeder&Platten99}. Moreover,

\begin{equation}
\label{F:MSE_error}
MSE_{N,h}[\phi_h]:=E\big[\big(\phi_h-E[\phi_h]\big)^2\big]\leq B_{\Xi}^2(E[\phi]) + V[\phi]/N, 
\end{equation}

where MSE stands for mean square error. Since
$E^2[|\phi_h-E[\phi]|]\leq E[(\phi_h-E[\phi_h])^2]$, $E[\phi_h]=u_h({\bf x}_0)$, and $E[\phi]=u({\bf x}_0)$, it holds 

\begin{eqnarray}
\label{F:Bound_of_error}
|u({\bf x}_0)-u_h({\bf x}_0)|\leq\sqrt{MSE(\phi_h)}=\sqrt{\big(|B_{\Xi}(E[\phi_h])|+\sqrt{V[\phi]/N}\big)^2
-2|B_{\Xi}^2(E[\phi_h])|V[\phi]/N}\leq \\
\leq |B_{\Xi}(E[\phi_h])|+\sqrt{V[\phi]/N}
\leq |B_{\Xi}(E[\phi_h])|+q\sqrt{V[\phi]/N},\nonumber
\end{eqnarray}

with probabilities approaching $P_q$ as $h\shortrightarrow 0^+$. Looking at (\ref{F:MSE_error}) or (\ref{F:Bound_of_error}), the two ways to speed up Monte Carlo simulations of (\ref{F:Dynkin}) are: using an integrator $\Xi$ with the highest possible $\delta$ for the BVP under consideration; and/or replacing the mean in (\ref{F:Statistical_error}) by another estimator of the expected value which has less variance, thus achieving the same statistical error with a smaller $N$.

\subsection{Variance reduction based on pathwise control variates}\label{SS:VR}
In this paper, we investigate a technique of variance reduction based on pathwise control variates, which is difficult to implement in other contexts, but suits the framework of PDD ideally. In order to discuss it, we adopt the same ansatz as in \cite[chapter 2]{Milstein&Tretyakov_book}. Consider the system of SDEs:
\begin{equation}
\label{F:VR_Milstein_sys}
\left\{
\begin{array}{ll}
d{\bf X}_t= [b({\bf X}_t)-\sigma({\bf X}_t){\bm \mu}({\bf X}_t)]dt + \sigma({\bf X}_t)d{\bf W}_t, & {\bf X}_0={\bf x}_0, \\
dY_t= c({\bf X}_t)Y_tdt + {\bm \mu}^T({\bf X}_t)Y_td{\bf W}_t, & Y_0=1,\\
dZ_t= -f({\bf X}_t)Y_tdt + {\bf F}^T({\bf X}_t)Y_td{\bf W}_t, & Z_0= 0,
\end{array}
\right. 
\end{equation}

where ${\bm \mu},{\bf F}:\Omega\shortrightarrow {\mathbb R}^D$ are smooth arbitrary fields, and let $({\bf X}_{\tau},Y_{\tau},Z_{\tau})$ be the evaluation of the processes $({\bf X}_t,Y_t,Z_t)$ at $t=\tau$. If ${\bm \mu}({\bf x})={\bf F}({\bf x})=0$, the system (\ref{F:VR_Milstein_sys}) is just another way of writing down the functional in (\ref{F:Dynkin}):
\begin{equation}
g({\bf X}_{\tau})Y_{\tau}+Z_{\tau}=\left\{
\begin{array}{ll}
\phi, & \textrm{  if  } {\bm \mu}({\bf x})={\bf F}({\bf x})=0,\\
{\tilde\phi}, & \textrm{otherwise.}
\end{array}
\right.
\end{equation}
It turns out--see \cite{Milstein&Tretyakov_book} for details--that 
\begin{equation}
E[g({\bf X}_{\tau})Y_{\tau}+Z_{\tau}]=E[{\tilde\phi}]=E[\phi]
\end{equation}
regardless of the arbitrary functions ${\bm \mu}$ and ${\bf F}$, but the variance
\begin{equation}
\label{F:Variance_VRSystem}
V[\phi]\neq V[{\tilde\phi}]= E\Big[\, \int_0^{\tau}Y_t^2\parallel \sigma^T\nabla u+u{\bm\mu}+{\bf F}\parallel_2^2\,dt \Big]
\end{equation}

{\em does} depend on them; and in fact, making the choice
\begin{equation}
\label{F:Combining_method}
\sigma^T\nabla u+u{\bm\mu}+{\bf F} =0
\end{equation}

then $V[{\tilde\phi}]=0$--meaning that one single realization would yield the exact pointwise solution $u({\bf x}_0)$ deterministically. Certainly this is a moot point, for in order to do so, the exact solution would be required in the first place; and the system (\ref{F:VR_Milstein_sys}) would still have to be integrated numerically--so that the variance of the discretized approximation would be small, but finite.

Let ${\tilde u}({\bf x})\approx u({\bf x})$. Choosing ${\bf F}=0$ in (\ref{F:Combining_method}),
\begin{equation}
\label{F:Exactmu}
{\bm \mu}= -\frac{1}{{\tilde u}}{\bm \sigma}^T \nabla {\tilde u}
\end{equation}
leads to the pathwise equivalent of importance sampling, while if ${\bm \mu}=0$,
\begin{equation}
\label{F:ExactF}
{\bf F}= -{\bm\sigma}^T \nabla {\tilde u}
\end{equation}
can be interpreted as the stochastic equivalent of control variates. Both importance sampling and control variates are well-known variance reduction techniques in statistics \cite{Glasserman2004}. 
There are several reasons why the second method is the more convenient:
\begin{itemize}
\item Control variates do not need to store ${\tilde u}({\bf x})$.
\item For inadequate ${\tilde u}$ and $\nabla {\tilde u}$, importance sampling may actually lead to {\em increased} variance. This is much less likely so with control variates (as it will be seen in a moment). 
\item From the point of view of numerical integration, introducing a non-zero $\bm \mu$ in (\ref{F:VR_Milstein_sys}) leads to a system of stochastic equations with multiplicative noise: $dY_t= c({\bf X_t})Y_tdt + {\bm \mu}^T({\bf X_t})Y_td{\bf W}_t$. Such SDEs can be unstable \cite{Kloeder&Platten99}. 
\item Finally, the trajectories are not affected by control variates, while if $\bm \mu\neq 0$, the 'particles' are actually 'pushed' around depending on $\nabla {\tilde u}$. The first fact facilitates the analysis of IterPDD.
\end{itemize}

Therefore, we will settle for (\ref{F:ExactF}) henceforth. 
In statistics, a control variate $\xi$ for $\phi$ is a variable which is drawn alongside $\phi$ and has a large correlation with it (i.e. close to $\pm 1$). Then, the variable
\begin{equation}
\label{F:Statistics_VR}
{\tilde\phi}(\gamma)= \phi - \gamma(\xi - E[\xi] ) 
\end{equation}
is such that $E[{\tilde\phi}]=E[\phi]$, but the variance is not. At the unique critical point $\gamma_*=cov[\phi,\xi]/V[\xi]$, $V[{\tilde\phi}]$ is minimal and yields a reduction of variance 
\begin{equation}
\label{F:Variance_Reduction}
\frac{V[{\tilde \phi}]}{V[\phi]} = 1-\rho^2[\phi,\xi] 
\end{equation}
(see \cite{Glasserman2004} for details), where $\rho$ is Pearson's correlation:
\begin{equation}
\label{F:Rho}
\rho[\phi,\xi]:= \frac{Cov[\phi,\xi]}{\sqrt{V[\phi]V[\xi]}}\qquad\textrm{(such that $\rho^2[\phi,\xi]\leq 1$).}
\end{equation}

By comparison of (\ref{F:Statistics_VR}) and (\ref{F:VR_Milstein_sys}), the pathwise control variate is Ito's integral 
\begin{equation}
\label{F:xiControlVariate}
\xi:= -\int_0^{\tau}Y_t{\bm \sigma}^T({\bf X}_t)\nabla {\tilde u}\, d{\bf W}_t. 
\end{equation} 

(Note that $k\xi$, with $k\neq 0$, yields the same $\rho[\phi,\xi]$ as $\xi$ and thus the same reduction of variance.) When ${\tilde u}=u$, we will refer to $\xi$ as the {\em exact control variate}. In that case, the correlation is perfect and the variance is zero (in the limit $h\shortrightarrow 0^+$). If ${\tilde u}\approx u$, $\xi$ is close to the minimizer $\gamma_*$ of (\ref{F:Statistics_VR}) and the equality in (\ref{F:Variance_Reduction}) still holds. But even if ${\tilde u}$ is so poor that (\ref{F:Variance_Reduction}) breaks down becoming a mere lower bound, any reasonable substitute of $u$ will still produce a $\xi$ with sufficient correlation as not to increase the variance. Summing up, the method of pathwise control variates is intrinsically robust--unlike pathwise importance sampling (\ref{F:Exactmu}). The control variate $\xi$ can be drawn by a quadrature of Ito's integral (\ref{F:xiControlVariate}), or alternatively, by enlarging system (\ref{F:VR_Milstein_sys}) with one extra equation:
\begin{equation}
d\xi_t= -Y_t{\bm \sigma}^T({\bf X}_t)\nabla {\tilde u}({\bf X}_t)\, d{\bf W}_t,\qquad \xi_0=0.
\end{equation}

Importantly, the same pathwise control variate can be used for parabolic BVPs and BCs involving derivatives. We close this section by adressing the choice of integrator $\Xi$ for solving (\ref{F:VR_Milstein_sys}) with ${\bf F}=-{\bm\sigma}^T \nabla {\tilde u}$ and ${\bm\mu}=0$. The simplest integrator for bounded SDEs is the Euler-Maruyama scheme plus a naive boundary test. This yields $\delta=1/2$, which is much poorer than the linear rate of weak convergence of the Euler-Maruyama integrator in free space \cite{Kloeder&Platten99}. However, if the solution, coefficients, and boundary of the BVP are smooth enough, then there are a variety of methods which manage to raise $\delta$ to one--see \cite{OurSurvey} for a recent review. In particular, the integrator of Gobet and Menozzi \cite{Gobet&Menozzi_2010} restores weak linearity of the Euler-Maruyama scheme by simply shrinking the domain (Algorithm \ref{A:Gobet_Menozzi}). (The reader is referred to that paper as to why the local shrinkage is that on line 3 below.)

\begin{algorithm}[h!]
\caption{Gobet and Menozzi's integrator with pathwise control variates.}
\begin{algorithmic}[1]
\STATE{{\bf Data:} $h>0$, ${\bf X}_0={\bf x}_0\in\Omega,{Y_0}=1\xi_0=0,Z_0=0$, $(d_0<0,{\bf N}_0)$}
\FOR{$k=0,1,2,\ldots$ until hitting the shrunken boundary}
	\IF{$d_k>-0.5826||\sigma_k{\bf N}_k||\sqrt{h}$}
		\STATE{Take one unbounded step according to:}
		\begin{eqnarray}
		\label{F:E_M}
		\left\{\begin{array}{l}
		{\bf X}_{k+1}= {\bf X}_k + h{\bf b}_k + \sqrt{h}\sigma_k{\cal N}_D, \\
		Y_{k+1}= Y_k + hY_kc_k,  \\
      \xi_{k+1}= \xi_k + hY_k{\bm\sigma}^T\nabla {\tilde u}({\bf X}_k),\\
	   Z_{k+1}= Z_k + hY_kf_k + \xi_{k+1}.
		\end{array}\right.
		\end{eqnarray}
		\STATE{Compute $(d_{k+1},{\bf N}_{k+1})$ according to the distance map of $\Omega$}
	\ELSE{}
		\STATE{Finish: $\tau_h=kh,{\bf X}_{\tau}=\Pi_{\partial\Omega}({\bf X}_k),\, \xi=\xi_k$ and $\phi_h= g({\bf X}_{\tau_h})Y_{\tau_h}+Z_{\tau_h}$.}
	\ENDIF
\ENDFOR
\end{algorithmic}
\label{A:Gobet_Menozzi}
\end{algorithm}
 
In Algorithm \ref{A:Gobet_Menozzi}, $d({\bf x})$ is the signed distance from ${\bf x}\in{\mathbb R}^D$ to $\partial\Omega$ with the convention that it is negative inside $\Omega$ (see \cite{Chip} for further details as to how to produce it); $\tau_h$ is the approximation to $\tau$; ${\cal N}_D$ is a $D-$dimensional standard normal distribution; $\Pi_{\partial\Omega}({\bf x})$ is the closest point on $\partial\Omega$ to ${\bf x}$, and ${\bf N}$ is the unit outward normal at $\Pi_{\partial\Omega}({\bf x})$ ($\Omega$ is supposed smooth enough that it exists).

\section{Cost and correlation of the control variates in PDD} \label{S:Theory}
\subsection{Complexity of a PDD simulation}\label{SS:Cost}

\begin{definition}
\label{D:Balanced_MC_simulation}
A balanced MC simulation with accuracy $a$ and confidence interval $P_q$ is such that: statistical error = $\frac{a}{2}$ = absolute value of signed bias.
\end{definition}

A balanced MC simulation is thus one guaranteed to attain a total $MSE_{N,h}(u_h(\bf x)_0)$ (\ref{F:MSE_error}) smaller than $a$ (with a probability $P_q$) with the least computational complexity (cost), because it takes the largest possible $h$ and the fewest possible realizations compatible with $a/2$ and $P_q$, namely

\begin{equation}
\label{F:N_y_h}
N=\frac{4q^2V[\phi]}{a^2}\qquad\textrm{  and  }\qquad h=\big( \frac{a}{2|\beta|} \big)^{1/\delta}.
\end{equation}

(We remark in passing that, while splitting the total error between bias and variance looks natural enough, it is definitely suboptimal \cite{Haji-Ali}.) 
The average number of steps before hitting the boundary is
\begin{equation}
\label{F:nu}
\nu= E[\tau]/h,
\end{equation} 

where mean first exit time, $E[\tau]$, is well defined for a given BVP and ${\bf x}_0$. On average, every trajectory makes $\nu$ visits to the integrator $\Xi$. 
For a balanced MC simulation and a tolerance $a$, the expected cost is, then:
\begin{equation}
\label{F:Pointwise_cost}
\textrm{nodal cost}({\bf x}_0)= N\nu= NE[\tau]/h= 4q^2E[\tau](2|\beta|)^{1/\delta}\frac{V[\phi]}{a^{2+1/\delta}}:=
K\frac{V[\phi]}{a^{2+1/\delta}}.
\end{equation}

Using control variates to reduce the variance, the trajectories remain the same, but the cost per timestep is increased by a factor $\kappa\gtrsim 1$. This is due to 
the extra cost of interpolating the control variates from a lookup table (a $D-$ dimensional array holding pointwise values) and evaluating the extra term ${\bf F}$ in (\ref{F:VR_Milstein_sys}). Inserting (\ref{F:Variance_Reduction}), the pointwise cost with control variates is

\begin{equation}
\label{F:VR_Pointwise_cost}
\textrm{nodal cost with control variates}({\bf x}_0)= \kappa K\frac{V[\phi+\xi](1-\rho^2[\phi,\xi])}{a^{2+1/\delta}}.
\end{equation}

The cost of 
a global PDD approximation $u_0({\bf x})$ with PlainPDD($a_0$) is then

\begin{equation}
\label{F:PlainPDD_cost}
\textrm{Cost of PlainPDD($a_0$)}= \Pi + \sum_{i=1}^n \frac{K_iV_i[\phi]}{a_0^{2+1/\delta}},
\end{equation}
 
$\Pi$ is the cost (independent of $a_0$) involved in solving all the subdomains using the deterministic subdomain solver once the interfacial values are available. Finally, 

\begin{equation}
\label{F:IterPDD_cost}
\textrm{Cost of IterPDD($a_0,a_1$)}= 2\Pi + {\tilde \Pi} + \frac{1}{a_0^{2+1/\delta}}\sum_{i=1}^n K_iV_i[\phi]
\Bigg( \big(\frac{a_0}{a_1}\big)^{2+1/\delta} + \kappa\big(1-\rho_i^2[\phi,\xi]\big) \Bigg).
\end{equation}

Above, ${\tilde \Pi}$ is the cost (independent of $a_0$ and $a_1$) of constructing and storing ${\bf F}=-{\bm \sigma}^T\nabla {\tilde u}_1(\bf x)$ with ${\tilde u}_1=$PlainPDD$(a_1)$, and $\rho_i[\phi,\xi]$ depends on $a_1$ only as long as there are no quadrature errors (ie as $h\shortrightarrow 0^+$). 

\subsection{The correlation $\rho[\phi,\xi(a)]$ in the limit $h\shortrightarrow 0^+$}\label{SS:Lemmas}

The results in this subsection are exact (within the assumptions) in the limit $h\shortrightarrow 0^+$; approximations are left to Section \ref{S:Approximations}. To the best of our knowledge, they are also new. Let us introduce the following notation:
\begin{itemize}
\item $\eta\sim P$, with $\eta\in{\mathbb R}^s$ means that $\eta$ is some realization of the distribution $P$. The notation $\eta=P({\bm\omegaup})$ denotes a specific realization, labeled with the 'chance variable' ${\bm\omegaup}\in{\mathbb R}^s$.  
Thus, $\eta=P({\bm\omegaup})$ is the same as $\eta\sim P$, but the former considers ${\bm\omegaup}$ fixed and treats $\eta$ as a scalar (or vector). 
\item ${\cal N}$ is an $s$-dimensional standard normal distribution, with $s=1$ (by default), or known from the context. For instance, ${\cal N}({\bm\omegaup})$ is the realization labeled by ${\bm\omegaup}\in{\mathbb R}^s$.
\item The notation $\eta(\cdot;{\bm \omegaup})$ means that the stochastic variable $\eta$ (which may depend on several parameters represented by the first argument) is constructed based on a PDD simulation with nodal statistical errors labeled with $\bm\omegaup=(\omega_1,\ldots,\omega_n)$, where $\omega_i,\,1\leq i \leq n$ labels the statistical error on node ${\bf x}_i$. (An arbitrary node is denoted by ${\bf x}_0$ and $\omega_0$.) 
\end{itemize}   

\begin{lemma}
\label{Th:Lema1}
Let $\phi$ and $\eta$ be stochastic variables, $\xi$ an exact control variate for $\phi$, and assume that $V[\phi],V[\xi]$ and $V[\eta]$ are finite. It holds:
\begin{enumerate}
\item $V[\xi]=V[\phi]$.
\item $Cov[\phi,\xi]=-V[\phi]\leq 0.$
\item $Cov[\xi,\eta]=-Cov[\phi,\eta]$.
\end{enumerate} 
\end{lemma}
{\em Proof.} We will make use of the Cauchy-Schwarz inequality $|Cov[a,b]|\leq \sqrt{V[a]V[b]}$ for $a$,$b$ with finite variances. Notice first that
$V[\phi+\xi]=0=V[\phi]+V[\xi]+Cov[\phi,\xi]$, so that $Cov[\phi,\xi]$ cannot be positive. For the first result, note that $|Cov[\phi,\xi]|=-Cov[\phi,\xi]$, and using the Cauchy-Schwarz inequality, $-2\sqrt{V[\phi]V[\xi]}-2Cov[\phi,\xi]\leq 0$. Summing this and $V[\phi]+V[\xi]+2Cov[\phi,\xi]=0$ yields
\begin{equation}
0\geq V[\phi] + V[\xi] - 2\sqrt{V[\phi]V[\xi]}= \big(\sqrt{V[\phi]} - \sqrt{V[\xi]} \,\big)^2 
\end{equation}  
which can only happen if $V[\phi]=V[\xi]$. Then, it is clear that $2Cov[\phi,\xi]=-2V[\phi]$. Finally, 
$|Cov[\phi+\xi,\eta]|\leq \sqrt{V[\phi+\xi]V[\eta]}=0$, so that 
$Cov[\phi+\xi,\eta]=Cov[\phi,\eta]+Cov[\xi,\eta]=0$ and the third result follows. 
$\Square$\newline

The central limit theorem behind the estimates for the statistical error in (\ref{F:Statistical_error}) assumes that Monte Carlo errors are normally distributed. Moreover, they are biased due to discretization. Lemma \ref{Th:Lema3} makes this assumption explicit.

\begin{lemma}
\label{Th:Lema3}
Let $u_h\sim u + \beta h^{\delta} + \sqrt{V[\phi]/N}{\cal N}$. Then $E[u_h-u]=\beta h^{\delta}$ and $E[(u_h-u)^2]=\beta^2h^{2\delta}+V[\phi]/N$.
\end{lemma}
{\em Proof.} The first moment is trivial. For the MSE, just notice that $E[(u_h-u)^2]= E[(\beta h^{\delta}+\sqrt{V[\phi]/N}{\cal N})^2]= \beta^2h^{2\delta} + (V[\phi]/N)E[{\cal N}^2] + \beta\sqrt{V[\phi]/N}E[{\cal N}]$. \Square\newline

The distribution of $u_h$ defined in Lemma \ref{Th:Lema3} results from combining the central limit theorem in (\ref{F:Statistical_error}) with the $MSE(\phi_h)$ in (\ref{F:MSE_error}). It is natural enough and further justified by the match of the first two moments of the error, but rigourous only asymptotically as $h\shortrightarrow 0^+$.

In order to track errors from the nodal values into the subdomains, one needs to know how they are propagated by the interfacial interpolators. Definition \ref{D:Linear_interpolator} introduces the relevant notation and properties. 

\begin{definition}
\label{D:Linear_interpolator}
Let $\lambda,\mu$ be real constants, $\{{\bf x}_1,\ldots,{\bf x}_r\}$ a set of $r>1$ distinct points in $\Gamma\subset{\mathbb R}^D$;
$u_1,\ldots,u_r$ and $v_1,\ldots,v_r$ two sets of scalars associated to the points in $\Gamma$, and $z:\Gamma\mapsto{\mathbb R}$ a function.
Let $R[z({\bf x}_i)\,|\,{\bf x}_i\in\Gamma]({\bf x}):\Gamma\mapsto{\mathbb R}$ be a smooth approximation to $z$ obtained by interpolation of $\{z({\bf x}_1),\ldots,z({\bf x}_r)\}$. We say that $R$ is a linear interpolator if 
\begin{equation}
R[\lambda u_i + \mu v_i\,|\,{\bf x}_i\in\Gamma]({\bf x})=
\lambda R[u_i\,|\,{\bf x}_i\in\Gamma]({\bf x}) + \mu R[v_i\,|\,{\bf x}_i\in\Gamma]({\bf x}), \qquad {\bf x}\in{\Gamma}.
\end{equation}
For some norm $||.||$ in $\Gamma$ we call the interpolation error
\begin{equation}
\epsilonup_z:=||R[z({\bf x}_i)\,|\,{\bf x}_i\in\Gamma]-z||.
\end{equation}   
\end{definition}

The most common interpolation schemes are linear in the sense of Definition \ref{D:Linear_interpolator}--for instance, RBF interpolation \cite{Fasshauer_libro}.

\begin{definition}
\label{D:Error_propagation_function}
The error propagation function $w_k({\bf x};{\bm\omegaup})$ for the elliptic BVP (\ref{F:EllipticBVP}) in the subdomain $\Omega_k\,(1\leq k\leq m)$ is defined as
\begin{equation}
\label{F:Error_propagation_function}
w_k({\bf x};{\bm\omegaup})= {\bar w}_k({\bf x}) + \frac{1}{q}{\tilde w}_k({\bf x};{\bm\omegaup})
\end{equation}
where ${\bar w}_k({\bf x})$ and ${\tilde w}_k({\bf x};{\bm\omegaup})$ are respectively the solution of the deterministic BVP and of the BVP with stochastic BCs
\begin{equation}
\label{F:w_bar}
\left\{
\begin{array}{lll}
L{\bar w}_k({\bf x})+c{\bar w}_k({\bf x})= 
L{\tilde w}_k({\bf x};{\bm\omegaup})+c{\tilde w}_k({\bf x};{\bm\omegaup})= 0,&\textrm{ if }&{\bf x}\in\Omega_k, \\
{\bar w}_k({\bf x})= {\tilde w}_k({\bf x};{\bm\omegaup})= 0,&\textrm{ if }&{\bf x}\in\partial\Omega_k\cap\partial\Omega, \\
\left.
\begin{array}{l}
{\bar w}_k({\bf x})= R[\,sign(\beta_i)\,|\,{\bf x}_i\in\Gamma_j^k\,]({\bf x}) \\
{\tilde w}_k({\bf x};{\bm\omegaup})= R[\,{\cal N}(\omega_i)\,|\,{\bf x}_i\in\Gamma_j^k\,]({\bf x})
\end{array}
\right\}
&\textrm{ if }&{\bf x}\in\Gamma_j^k,\,1\leq j \leq m_k.
\end{array}
\right.
\end{equation}
\end{definition}


\begin{definition}
\label{D:definition_of_psi}
Let $\tau$, $\sigma$, ${\bf X}_t$, $Y_t$, and ${\bf W}_t$ be the same as in (\ref{F:xiControlVariate}). For a fixed ${\bm\omega}\in{\mathbb R}^n$,
\begin{equation}
\label{F:definition_of_psi}
\psi(\bm\omega):= -\int_{t=0}^{\tau}{\sigma}^T({\bf X}_t)\Bigg(\, \oplus_{k=1}^m \nabla\big( {\bar w}_k({\bf X}_t) + \frac{1}{q}{\tilde w}_k({\bf X}_t;{\bm\omega})\big)\,\Bigg)Y_td{\bf W}_t,
\end{equation}
\end{definition}

Thanks to the smoothness of the iterfacial interpolator, gradients inside $\Omega_k$ are well defined, and thus the integral in (\ref{F:definition_of_psi}) is also well defined regardless of the continuity of $\oplus_{k=1}^m\nabla w_k({\bf x};{\bm\omega})$ across the interfaces. We are now prepared to state the main theoretical result.

\begin{lemma}
\label{Th:Teorema}
Assume an elliptic BVP with Dirichlet BCs like (\ref{F:EllipticBVP}) with exact solution $u({\bf x})$, and let ${\tilde u}({\bf x},a;{\boldsymbol \omega})$ be a PDD simulation of it in the limit $h\shortrightarrow 0^+$, with accuracy $a>0$, nodal statistical errors labeled by ${\bm\omegaup}$, balanced MC simulations, and smooth linear interpolator $R$. Let $\xi(a;{\boldsymbol \omega})$ be the control variate at a given interfacial node ${\bf x_0}\in\{{\bf x}_1,\ldots,{\bf x}_n\}$ constructed from ${\tilde u}({\bf x},a;{\boldsymbol \omega})$ according to (\ref{F:xiControlVariate}). Then
\begin{equation}
\label{F:V[xi(a,w)]}
V[\xi(a;{\bm\omegaup})] = V[\phi]-Cov[\phi,\psi(\bm\omegaup)]\,a+V[\psi(\bm\omegaup)]\big(\frac{a}{2}\big)^2
+{\cal O}\Bigg(\sum_{k=1}^m\sum_{j=1}^{m_k}\epsilonup_u|_{\Gamma_j^k}\Bigg)
+{\cal O}\Bigg(\sum_{k=1}^m\epsilonup_{\Omega_k}(a)\Bigg)
\end{equation}
and
\begin{equation}
\label{F:rho(a,w)}
\rho[\phi,\xi(a;{\bm\omegaup})] = 
-\sqrt{\frac{V[\phi]}{V[\xi(a;{\bm\omegaup})]}} + \frac{Cov[\phi,\psi({\bm\omegaup})]}{\sqrt{V[\phi]V[\xi(a;{\boldsymbol\omega})]}}\frac{a}{2}
+{\cal O}\Bigg(\sum_{k=1}^m\sum_{j=1}^{m_k}\epsilonup_u|_{\Gamma_j^k}\Bigg)
+{\cal O}\Bigg(\sum_{k=1}^m\epsilonup_{\Omega_k}(a)\Bigg),
\end{equation}
where $\psi(\boldsymbol\omega)$ is the auxiliar variate defined by (\ref{F:definition_of_psi}), $\epsilonup_u|_{\Gamma}$ the error of interpolating $u$ along an interface $\Gamma$, and $\epsilonup_{\Omega}$ is the error of the PDD subdomain solver.
\end{lemma}

{\em Proof.} Let us consider the integral representation of the solution of (\ref{F:EllipticBVP})
\begin{equation}
\label{F:Green_representation}
u({\bf x})= \int_{\Omega}G({\bf x},{\bf y})f({\bf y})d^D{\bf y} +
\int_{\partial\Omega}\frac{\partial G({\bf x},{\bf y})}{\partial {\bf N}}g({\bf y})d^{D-1}{\bf y},
\end{equation}
where $\partial/\partial{\bf N}={\bf N}\cdot\nabla$ and $G({\bf x},{\bf y})$ is Green's function, defined as the solution of
\begin{equation}
\label{F:Green_G}
\left\{
\begin{array}{lcr}
LG({\bf x},{\bf y})+cG({\bf x},{\bf y})=\delta({\bf x}-{\bf y}) &\textrm{ if }& {\bf x}\in\Omega,\\
G({\bf x},{\bf y})=0 &\textrm{ if }& {\bf x}\in\partial\Omega.
\end{array}
\right.
\end{equation} 
Under the adequate smoothness requirements on $L$ and $\partial\Omega$ the solution $G({\bf x},{\bf y})$ to the homogeneous BVP (\ref{F:Green_G}) exists and is unique, which ensures the validity of (\ref{F:Green_representation}) \cite{Bers1964}. Consider the solution restricted to subdomain $\Omega_k$. Since $G({\bf x},{\bf y})$ does not depend on the boundary data, $u|_{\Omega_k}({\bf x})$ can also be represented as
\begin{equation}
\label{F:Greens_representation_of_uk}
u|_{\Omega_k}({\bf x})=
\int_{\Omega_k}G({\bf x},{\bf y})f({\bf y})d^D{\bf y} +
\int_{\partial\Omega_k\cap\partial\Omega}\frac{\partial G({\bf x},{\bf y})}{\partial {\bf N}}g({\bf y})d^{D-1}{\bf y} +
\sum_{j=1}^{m_k}\int_{\Gamma_j^k}\frac{\partial G({\bf x},{\bf y})}{\partial {\bf N}}u({\bf y})d^{D-1}{\bf y}.
\end{equation}
On the other hand, the PDD subdomain approximation is affected by the error of the subdomain solver, $\epsilonup_{\Omega_k}(a)$, which depends on $a$ due to the BC along the interfaces: 
\begin{eqnarray}
\label{F:Greens_representation_of_vk}
{\tilde u}|_{\Omega_k}({\bf x},a;{\bm \omegaup})=
\int_{\Omega_k}G({\bf x},{\bf y})f({\bf y})d^D{\bf y} +
\int_{\partial\Omega_k\cap\partial\Omega}\frac{\partial G({\bf x},{\bf y})}{\partial {\bf N}}g({\bf y})d^{D-1}{\bf y} +\\\nonumber
+\sum_{j=1}^{m_k}\int_{\Gamma_j^k}\frac{\partial G({\bf x},{\bf y})}{\partial {\bf N}}{\tilde u}({\bf y},a;{\bm \omegaup})d^{D-1}{\bf y}+\epsilonup_{\Omega_k}(a).
\end{eqnarray}
In the limit $h\shortrightarrow 0^+$, the nodal values are given by the distribution in Lemma \ref{Th:Lema3}. Running a PDD simulation with balanced MC simulations, tolerance $a$ and probability $P_q$, is then equivalent to fixing ${\bm\omegaup}$ and taking
\begin{equation}
\label{F:Fake_nodal_values}
{\tilde u}({\bf x}_0,a;{\bm\omegaup})= u({\bf x}_0)+\frac{a}{2}sign\big(\beta({\bf x}_0)\big) + \frac{a}{2q}{\cal N}(\omega_0).
\end{equation}
(Note that the statistical error, and hence ${\bm\omegaup}$, are only known after running the PDD simulation.) The Dirichlet BC condition on an interface $\Gamma$ is then
\begin{equation}
\label{F:v_restricted_to_Gamma}
{\tilde u}|_{\Gamma}({\bf x},a;{\bm\omegaup})= R[{\tilde u}_i\,|\,{{\bf x}_i\in\Gamma}]({\bf x})=
u|_{\Gamma}({\bf x})+\frac{a}{2}R\big[\,sign\big(\beta({\bf x}_i)\big)+\frac{1}{q}{\cal N}(\omega_i)\,|\,{{\bf x}_i\in\Gamma}\,\big]({\bf x}) + {\cal O}(\epsilonup_u|_{\Gamma}).
\end{equation}
Inserting (\ref{F:v_restricted_to_Gamma}) into (\ref{F:Greens_representation_of_vk}):

\begin{eqnarray}
{\tilde u}|_{\Omega_k}({\bf x},a;{\bm\omegaup})= u|_{\Omega_k}({\bf x}) + 
\frac{a}{2}\sum_{j=1}^{m_k}\int_{\Gamma_j^k}\frac{\partial G({\bf x},{\bf y})}{\partial{\bf N}}R\big[\,sign\big(\beta({\bf x}_i)\big)+\frac{1}{q}{\cal N}(\omega_i)\,|\,{\bf x}_i\in\Gamma_j^k\,\big]({\bf y})d^{D-1}{\bf y}+\\\nonumber  
+ {\cal O}\Bigg(\sum_{j=1}^{m_k}\epsilonup_u|_{\Gamma_j^k}\Bigg) 
+\epsilonup_{\Omega_k}(a).
\end{eqnarray}

Since $R$ is a smooth interpolator, gradients are well defined inside $\Omega_k$ and along the interfaces, but not necessarily across them. In order to circumvent this issue we take the direct sum of subdomain gradients, 
\begin{equation}
\oplus_{k=1}^m \nabla {\tilde u}({\bf x},a;{\bm\omega})= \oplus_{k=1}^m\nabla u({\bf x})|_{\Omega_k} + \frac{a}{2}{\bf V}({\bf x};{\bm\omega}) 
+{\cal O}\Bigg(\sum_{k=1}^m\sum_{j=1}^{m_k}\epsilonup_u|_{\Gamma_j^k}\Bigg)
+{\cal O}\big(\epsilonup_{\Omega_k}(a)\big),
\end{equation}

where
\begin{equation}
\label{F:aux1}
{\bf V}({\bf x};\bm\omega)= \oplus_{k=1}^m \nabla\Bigg( \,\sum_{j=1}^{m_k}\int_{\Gamma_j^k}\frac{\partial G({\bf x},{\bf y})}{\partial {\bf N}}R\big[\,sign\big(\beta({\bf x}_i)\big)+\frac{1}{q}{\cal N}(\omega_i)\,|\,{\bf x}_i\in\Gamma_j^k\,\big]({\bf y})d^{D-1}{\bf y}\Bigg).
\end{equation}

Note that the quantity in parentheses in (\ref{F:aux1}) is the Green representation 
of the BVP with solution $w_k({\bf x};{\bm\omegaup})={\bar w}_k({\bf x})+\frac{1}{q}{\tilde w}_k({\bf x};{\bm\omegaup})$ (\ref{F:Error_propagation_function}), by linearity. Now, according to definitions 
(\ref{F:xiControlVariate}) and (\ref{F:definition_of_psi}) for node ${\bf x}_0$,

\begin{equation}
\label{F:xi(a)}
\xi(a;{\bm\omega})=\xi + \frac{a}{2}\psi(\bm\omega) + {\cal O}\Bigg(\sum_{k=1}^m\sum_{j=1}^{m_k}\epsilonup_u|_{\Gamma_j^k}\Bigg)+ \sum_{k=1}^m\epsilonup_{\Omega_k}(a).
\end{equation}

In (\ref{F:xi(a)}), we have used the fact that $\oplus_{k=1}^m\nabla u|_{\Omega_k}=\nabla u$ except on the interfaces. Since the interfaces are smooth and the random paths are not, the trajectories cross the interfaces at isolated points which make no contribution to the integral, so that $\displaystyle -\int_{t=0}^{\tau}{\bm\sigma}^T\big(\oplus_{k=1}^m\nabla u|_{\Omega_k}\big)Y_td{\bm W}_t=\xi$. Taking the variance of (\ref{F:xi(a)}) and using Lemma \ref{Th:Lema1} yields (\ref{F:V[xi(a,w)]}). Moreover,

\begin{equation}
Cov[\phi,\xi(a;\boldsymbol\omega)]= Cov[\phi,\xi]
+\frac{a}{2}\,Cov[\phi,\psi({\boldsymbol\omega)}]
+{\cal O}\Bigg(\sum_{k=1}^m\sum_{j=1}^{m_k}\epsilonup_u|_{\Gamma_j^k}\Bigg) + {\cal O}\Bigg(\sum_{k=1}^m\epsilonup_{\Omega_k}(a)\Bigg).
\end{equation}\newline


Finally, equation (\ref{F:rho(a,w)}) follows from $\rho[\phi,\xi(a;{\boldsymbol\omega})]=Cov[\phi,\xi(a;\boldsymbol\omega)]/\sqrt{V[\phi]V[\xi(a;{\boldsymbol\omega})]}$ recalling that the correlation with the exact control variate is $-1$ by Lemma \ref{Th:Lema1}. $\Box$\\


The point of Lemma \ref{Th:Teorema} is to predict the correlation between the score $\phi$ and an approximate control variate $\xi(a;{\bm \omegaup})$ without actually running a PDD simulation to produce the latter--but rather "simulating" it. In exchange, the variable $\psi({\bm\omegaup})$ must be computed, but it can be constructed on a subdomain-per-subdomain basis, i.e. in a fully parallellizable way. The drawback is that the formulas derived so far are only rigourous at $h=0$, and that they depend on quite a few problem-dependent constants--many of them node-dependent as well. These difficulties will be 
addressed in Section \ref{S:Approximations}. But before that, we shall examine the issues of global error and stability of PDD, and introduce some more results which will later be useful in assessing the necessary simplifications.

\section{Global error and stability of PDD} \label{S:GlobalError}

The main tool is Theorem 3.7 in Gilbarg and Trudinger \cite{G&T_book}:
\begin{theorem}
\label{Th:GandT}
Let $v$ be the solution of an elliptic BVP like (\ref{F:EllipticBVP}) with $c\leq 0$, such that $v$ is continuous on $\Omega\backslash\partial\Omega$ and twice differentiable on $\partial\Omega$. Then
\begin{equation}
\label{F:GandT}
\sup_{{\bf x}\in\Omega} |v|\leq \sup_{{\bf x}\in\partial\Omega} |g| + Q\sup_{{\bf x}\in\Omega} \frac{|f|}{\lambda} 
\end{equation}
where $Q$ is a constant depending only on $diam(\Omega)$ and $\sup_{\Omega}||{\bf b}||_2/\lambda$. In particular, if $\Omega$ lies between two parallel planes a distance $d$ apart, then (\ref{F:GandT}) is satisfied with $Q=\exp{[(\sup_{\Omega}||{\bf b}||_2/\lambda)/d]}-1$.
\end{theorem}
Recall that $\lambda=\inf_{\bf x\in\Omega}\lambda_{min}[A({\bf x})]$, and $diam(\Omega)$ is the largest distance between two points in $\Omega$. 
It is convenient to introduce the following parameter:

\begin{definition}
\label{D:Overshoot}
Let $R$ be a linear interpolator and $\Gamma$ a subdomain interface with nodes ${\bf x}_i\in\Gamma,\,1\leq i\leq p$. Let $z_1,\ldots,z_p$ be $p$ scalars such that $|z_i|= 1$. The overshoot constant of $R$ with respect to the discretization ${\bf x}_1,\ldots,{\bf x}_p$ of $\Gamma$ is defined as  
\begin{equation}
\label{F:Overshoot}
\gamma_R^{\Gamma}:= \sup\limits_{z_1,\ldots,z_p}\sup\limits_{{\bf x}\in\Gamma} \big|R[z_i|{\bf x}_i\in\Gamma]({\bf x})\big|.
\end{equation}
Analogously, let $\gamma_R^{\Omega_k}:=\sup_{\Gamma\in\Omega_k} \gamma_R^{\Gamma}$ and $\gamma_R:=\sup_{1\leq k\leq m}\gamma_R^{\Omega_k}$.
\end{definition}

The overshoot constant measures the excess of the reconstructed function over any of the interpolated values which are $1$ in absolute value. For piecewise interpolation, $\gamma_R=1$, but for more accurate interpolators, $\gamma_R>1$ due to the Runge and Gibbs phenomena \cite{Fasshauer_libro}--see Figure \ref{I:Figura3} (right) for illustration.  

We are interested in bounding the largest PDD error throughout $\Omega$. Assuming as always balanced MC simulations, and neglecting the error of the subdomain solver, the PDD error in subdomain $\Omega_k$ obeys
\begin{equation}
\label{F:Error_PDE}
\left\{
\begin{array}{ll}
L({\tilde u}-u)|_{\Omega_k}+c({\tilde u}-u)|_{\Omega_k}=0& \textrm{if ${\bf x}\in\Omega_k$,}\\
({\tilde u}-u)|_{\Omega_k}= 0& \textrm{if ${\bf x}\in\partial\Omega_k\cap\partial\Omega$,}\\
({\tilde u}-u)|_{\Omega_k}= \frac{a}{2}R[sign(\beta_i)+\frac{1}{q}{\cal N}(\omegaup_i)\,|\,{\bf x}_i\in\Gamma_j^k]({\bf x})& \textrm{if ${\bf x}\in\Gamma_j^k,\,1\leq j\leq m_k.$} 
\end{array}
\right.
\end{equation}

Therefore, $({\tilde u}-u)|_{\Omega_k}= ({\tilde u}_k-u|_{\Omega_k})= \frac{a}{2}w_k({\bf x};{\bm\omegaup})$--hence the name of error propagation function. Applying Theorem \ref{Th:GandT} to (\ref{F:Error_PDE}) and by linearity of $R$,


\begin{equation}
\big| {\tilde u}_k-u|_{\Omega_k}\big| \leq 
\frac{aQ_k}{2}
\sup_{{\bf x}\in\partial\Omega_k\backslash\partial\Omega}\Big|\oplus_{\Gamma_j^k\in\partial\Omega_k} R[sign(\beta_i)+\frac{1}{q}{\cal N}(\omegaup_i)\,|\,{\bf x}_i\in\Gamma_j^k]({\bf x})\Big|
\leq
\frac{aQ_k\gamma_R^{\Omega_k}}{2}\Big( 1 + 
\frac{1}{q}\sup_{{\bf x}_i\in\Omega_k}|{\cal N}(\omegaup_i)|\Big),
\end{equation}

where $Q_k$ depends on the size and shape of $\Omega_k$ and $L|_{\Omega_k}$. 
Thanks to the symmetry of the standard normal distribution around zero, $\sup_{1\leq j\leq s}{|\cal N}(\omegaup_j)|=2\sup_{1\leq j\leq s}{\cal N}(\omegaup_j)$, so that the {\em global} PDD error (neglecting the error of the subdomain solver) can be bounded by

\begin{equation}
\label{F:Global_PDD_error_bound}
|{\tilde u}-u| \leq \sup_{1\leq k \leq m} \,
\frac{a\gamma_R^{\Omega_k}Q_k}{2}\Big( 1 + \frac{2}{q}\sup_{{\bf x}_i\in\partial\Omega_k}{\cal N}(\omegaup_i)\Big).
\end{equation}

Equation (\ref{F:Global_PDD_error_bound}) reflects that the global PDD error depends on PDE coefficients in (\ref{F:EllipticBVP}), on the shape and size $|\Omega_k|$ of the subdomains (typically, $|\Omega|/m$), and on the number of nodes per subdomain (typically, around $n/m$), rather than on $|\Omega|$ and $n$. Therefore, PDD is intrinsically stable provided that the number, shape and discretization of the subdomains are so chosen that the subdomain errors are controlled. 
 
We will now derive the nodal tolerance $a_0(\epsilonup)$ required to enforce a set global PDD error tolerance $\epsilonup>0$. 
Recall that $n_k$ is the total number of interfacial nodes sitting on $\partial\Omega_k$. The distribution $\sup_{{\bf x}_i\in\Omega_k}{\cal N}(\omega_i)=\sup_{1\leq j\leq n_k}{\cal N}(\omega_j)$ is an {\em extreme value} distribution. The value $n_k$ for which its maximum is less than $x$ with probability $P_q$ is to be extracted from
\begin{equation}
x=CDF[\sup_{1\leq j\leq n_k}{\cal N}(\omega_j)]^{-1}(P_q)
\textrm{  (such that }Pr[\sup_{1\leq j\leq n_k} {\cal N}(\omega_j)< x ]= P_q.)
\end{equation}
where $CDF[\cdot]$ and $CDF[\cdot]^{-1}$ stand for the cumulative distribution function (CDF) of a given distribution and its inverse, respectively. Let
\begin{equation}
k_{max}= \arg\max\limits_{1\leq k\leq m} \frac{a\gamma_R^{\Omega_k}Q_k}{2}\Big( 1 + \frac{2}{q}CDF[\sup_{{\bf x}_i\in\partial\Omega_k}{\cal N}(\omegaup_i)]^{-1}(P_q)\Big) 
\end{equation}
and define
\begin{equation}
Q_{max}=Q_{k_{max}}\gamma_R^{\Omega_{k_{max}}}/\gamma_R,\qquad 
s=n_{k_{max}}, \qquad
S_s= \sup_{1\leq j \leq s}{\cal N}(\omegaup_j)\,\,(i.i.d.).
\end{equation} 

Since the $s$ standard normals in $S_s$ are i.i.d., the CDF is

\begin{eqnarray}
\label{F:CDF_of_Ss}
CDF[S_s](x)= \int_{-\infty}^x S_s(t) dt=
Pr[\big({\cal N}_1\leq x\big)\cap \ldots \cap \big({\cal N}_s\leq x\big)]=
\big( CDF[{\cal N}](x) \big)^s.
\end{eqnarray}   

For the standard normal distribution
\begin{equation}
CDF[{\cal N}](x)= \frac{1}{2\pi}\int_{-\infty}^xe^{-t^2}dt=
\frac{1}{2}\big( 1+\textrm{erf}(x/\sqrt{2}) \big),
\end{equation}

where erf$(x)= (2/\sqrt{\pi})\int_{0}^x e^{-t^2}dt$ is the error function. 
The nodal target tolerance $a_0$ can then be related to the global PDD error tolerance by

\begin{equation}
\label{F:a(eps)}
\epsilonup= \frac{a_0\gamma_RQ_{max}}{2} \Big( 1 + \frac{2\sqrt{2}}{q}
\textrm{erf}^{-1}(2P_q^{1/s}-1)\Big). 
\end{equation}

\begin{figure}[h]
\centerline{\includegraphics[width=1.\columnwidth]{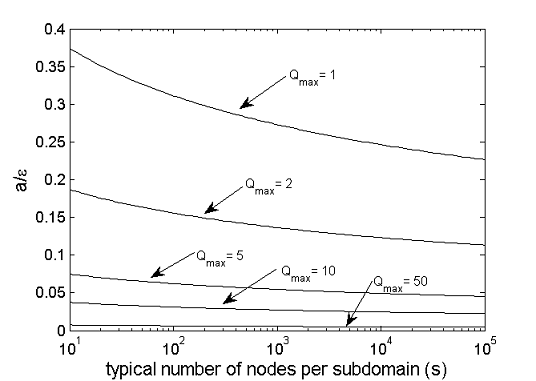}}
\caption{The ratio $a_0/\epsilonup$ from (\ref{F:a(eps)}) as a function of $s$ for several values of $Q_{max}$ and $q=2,\gamma_R=1.5$ (notice the different scales).}
\label{I:Figura2}
\end{figure}

As Figure \ref{I:Figura2} shows, $a_0/\epsilonup=a_0(\epsilonup=1)$ depends very mildly on the typical number of nodes per subdomain. Anyway, the bound (\ref{F:a(eps)}) will be a large overestimation in many cases, for the effect of the statistical errors decays fast away from the interfaces. If the moments of $a(\epsilonup)$ were required, a useful fact is that as $s$ grows, $S_s$ tends to the Gumbel distribution    
\begin{equation}
\lim_{s\shortrightarrow \infty} CDF[\frac{S_s-l_s}{b_s}](x)= G\big(\frac{x-l_s}{b_s}\big):= \exp[-e^{-x}]
\end{equation}
with location and scaling parameters $l_s=-CDF[{\cal N}]^{-1}(1/s)$ and $b_s= 1/l_s$.

\section{Approximations leading to a practical algorithm} \label{S:Approximations}

In order to construct an implementable multigrid-like IterPDD algorithm, several approximations are needed, listed as heuristics ${\bf H1}$ through ${\bf H5}$ below.\newline

\subsection*{H1: $a_0$ and $a>a_0$ are small enough}\label{SS:H1} 
This assumption is used on several occasions in {\bf H2}-{\bf H5}.\newline 

\subsection*{H2: Interpolation and subdomain-solver errors will be dropped}\label{SS:H2}

Interfacial interpolation errors in (\ref{F:V[xi(a,w)]}) and (\ref{F:rho(a,w)}) will be dropped. 
Then, inserting (\ref{F:V[xi(a,w)]}) into (\ref{F:rho(a,w)}) yields, after some manipulation, 
\begin{equation}
\label{F:aux2}
1-\rho^2[\phi,\xi(a;{\bm\omegaup})] \gtrsim
\frac{\frac{V[\psi({\bm\omegaup})]}{V[\phi]}\Big(1-\rho^2[\phi,\psi(\bm\omegaup)]\Big)(a/2)^2}
{1 -\sqrt{\frac{V[\psi(\bm\omegaup)]}{V[\phi]}}\rho[\phi,\psi(\bm\omegaup)]a
+\frac{V[\psi(\bm\omegaup)]}{V[\phi]}(a/2)^2},
\end{equation}

where the sign $\gtrsim$ has replaced the equality to make up for dropping the interpolation and subdomain-solver errors. Next, note that
\begin{equation}
\frac{V[\phi+\xi(a;{\bm\omegaup})]}{V[\phi]} \gtrsim 1-\rho^2[\phi,\xi(a;{\bm\omegaup})],
\end{equation}

since (\ref{F:Variance_Reduction}) only strictly holds for the minimizer of (\ref{F:Statistics_VR}), and with $a>0$, the PDD solution ${\tilde u}$ yields a control variate off the minimizer, regardless of $h$. Note also that the denominator in (\ref{F:aux2}) is positive, since $|\rho[\phi,\psi({\bm\omegaup})]|\leq 1$:

\begin{equation}
1 -\sqrt{\frac{V[\psi(\bm\omegaup)]}{V[\phi]}}\rho[\phi,\psi(\bm\omegaup)]a
+\frac{V[\psi(\bm\omegaup)]}{V[\phi]}(a/2)^2 \geq \Bigg(1-\sqrt{\frac{V[\psi({\bm\omegaup})]}{V[\phi]}}(a/2)\Bigg)^2.
\end{equation}


For small $a$ (see {\bf H1}), the denominator in (\ref{F:aux2}) can be dropped. More precisely, since the MC simulations are balanced, 
\begin{equation}
\sqrt{\frac{V[\psi({\bm\omegaup})]}{V[\phi]}}\frac{a}{2}=q\sqrt{\frac{V[\psi({\bm\omegaup})]}{N}},
\end{equation}

which is negligible if $V[\psi({\bm\omegaup})]<<N$--recall that this $N>>1$ is meant without variance reduction, and $V[\psi({\bm\omegaup})]={\cal O}(V[\bar\psi])$ is bounded by (\ref{F:V[barpsi]}) in {\bf H4} below. Assuming this and putting all together, it holds

\begin{equation}
\label{F:The_key}
\frac{V[\phi+\xi(a;{\bm\omegaup})]}{V[\phi]}\gtrsim \frac{V[{\psi({\bm\omegaup})}]a^2}{4V[\phi]}\Big( 1-\rho^2[\phi,{\psi({\bm\omegaup})}]\Big).
\end{equation}

The importance of (\ref{F:The_key}) is that the two factors affecting IterPDD--namely $a$ and the random statistical errors on the nodes--have been separated. Moreover, the latter has been expressed in terms of  $V[\psi({\bm\omegaup})]$ and $\rho[\phi,\psi({\bm\omegaup})]$. Also, 

\begin{itemize}
\item If the score $\phi$ and the auxiliar variate $\psi({\bm\omegaup})$ were perfectly correlated (i.e. $\rho^2[\phi,\psi(a,{\bm\omegaup})]=1$), then $V=[\phi+\psi(a,{\bm\omegaup})]=0$, meaning that $\psi(a,{\bm\omegaup})=k\xi$. Since this can only happen if $f=g=0$, the solution would be $u=0$.

\item The opposite limit, $\rho[\phi,\psi(a,{\bm\omegaup})]=0$, yields asymptotically
\begin{equation}
\lim_{a\shortrightarrow 0^+}\frac{V[\xi(a;{\bm\omegaup})]}{V[\phi]}= \frac{V[\psi({\bm\omegaup})]}{V[\phi]}(a/2)^2.
\end{equation}
\end{itemize}


\subsection*{H3: Small variance with respect to $\bm\omegaup$}\label{SS:H3}

The propagation of nodal statistical errors (labeled with ${\bm\omegaup}$) onto the subdomain solutions in critical in PDD. In Section \ref{S:GlobalError}, it was shown that the effect of ${\bm\omegaup}$ on the PDD aggregate error can be controlled on a subdomain-per-subdomain basis. In IterPDD($a_0,a_1$), there are two further aspects to ${\bm\omegaup}$. First, how much the variance reduction produced by ${\xi(a_1;{\bm\omegaup})}$ depends on the chance variable. Second, how it affects the decrease predicted by (\ref{F:The_key})--for $a_1$ will be determined based on that formula. In particular, the nodal simulations will in turn be "simulated" themselves by randomly drawing a chance variable--say ${\bm\omegaup'}$--and computing $V[\psi({\bm\omegaup'})]$ and $\rho[\phi,\psi({\bm\omegaup'})]$.

\begin{lemma}
\label{Th:E[w]}
Let $E_{\bm\omegaup}[\cdot]$ be the expected value relative to $\bm\omegaup$ and assume that interpolation errors are negligible. Then, $E_{\bm\omegaup}[{\tilde w}_k({\bf x};{\bm\omegaup})]=0$ ($1\leq k \leq m$).
\end{lemma} 

{\em Proof.} The Green function representation of ${\tilde w}_k$ is
\begin{equation}
{\tilde w}_k({\bf x};{\bm\omegaup})= \sum_{j=1}^{m_k}\int_{\Gamma_j^k}
\frac{\partial G({\bf x},{\bf y})|_{\Omega_k}}{\partial {\bf N}}R[{\cal N}({\omegaup_i})\,|\,{\bf x}_i\in\Gamma_j^k]({\bf y})d^{D-1}{\bf y},
\end{equation} 

where $G({\bf x},{\bf y})$ is determined by (\ref{F:Green_G}) and is deterministic. By linearity of the integral, the interpolator, and the expected value,
\begin{eqnarray}
E_{\bm\omegaup}[{\tilde w}_k({\bf x},{\bm\omegaup})]= \sum_{j=1}^{m_k}\int_{\Gamma_j^k}
\frac{\partial G({\bf x},{\bf y})|_{\Omega_k}}{\partial {\bf N}}R\Big[E_{\bm\omegaup}[{\cal N}({\omegaup_i})]\,|\,{\bf x}_i\in\Gamma_j^k\Big]({\bf y})d^{D-1}{\bf y}= \\\nonumber
\sum_{j=1}^{m_k}\int_{\Gamma_j^k}
\frac{\partial G({\bf x},{\bf y})|_{\Omega_k}}{\partial {\bf N}}
\big( 0 + \epsilon_0|_{\Gamma_j^k} \big)d^{D-1}{\bf y}, 
\end{eqnarray}
where $\epsilon_0|_{\Gamma_j^k}$ is the error in the reconstruction of the constant function $z=0$ on the interface $\Gamma_j^k$ with $R$, and which is zero by hypothesis. $\square$

As a consequence of Lemma \ref{Th:E[w]}, $ E_{\bm\omegaup}[\nabla{\tilde w}_k]=\nabla E_{\bm\omegaup}[{\tilde w}_k]=0$, and thus
\begin{equation}
\label{F:Bar_psi}
{\bar \psi}:= E_{\bm\omegaup}[\psi({\bm\omegaup})]= \int_0^{\tau}Y_t\sigma^T\sum\limits_{k=1}^m\nabla\big({\bar w}_k+\frac{1}{q}E_{\bm\omegaup}[{\tilde w}_k]\big)d{\bf W}_t= 
\int_0^{\tau}Y_t\sigma^T\sum\limits_{k=1}^m\nabla{\bar w}_kd{\bf W}_t.
\end{equation}


The variance can be calculated by Ito's isometry:

\begin{equation}
\label{F:Isometry}
V[\psi({\bm\omegaup})]= E\Big[\,\int_0^{\tau} Y_t^2({\bf X}_t)||\sigma^T({\bf X}_t)\sum\limits_{k=1}^m\nabla w_k({\bf X}_t,{\bm\omegaup})||^2_2\,dt\,\Big],
\end{equation}

where $Y_t= \exp{\int_0^{t}c({\bf X}_s)\,ds}$. 
By the triangular inequality and the inequalities

\begin{equation}
\label{F:Norm_matrixvector}
||\sigma^T\nabla {\tilde w}_k||_2\leq ||\sigma^T||_2\cdot||\nabla {\tilde w}_k||_2
\end{equation}

and

\begin{equation}
\label{F:Norm_of_matrix}
||\sigma^T({\bf x})||_2= \sqrt{\lambda_{max}\big(\sigma({\bf x})\sigma^T(\bf x)\big)}=\sqrt{2\lambda_{max}\big(A({\bf x})\big)}\leq \sqrt{2\Lambda},\textrm{ if ${\bf x}\in\Omega$,}
\end{equation}

one has

\begin{eqnarray}
\label{F:E[V[psi(w)]]}
E_{\bm\omegaup}\big[\, V[\psi({\bm\omegaup})] \,\big]\leq
V[{\bar \psi}] +
\frac{2\Lambda\,E\big[\int_0^{\tau}Y_t^2({\bf X}_t)dt\,\big]}{q^2}
\,E_{\bm\omegaup}\Big[\,\big(\sup_{{\bf x}\in\Omega}||\sum\limits_{k=1}^m{\nabla\tilde w}_k({\bf x};{\bm\omegaup})||_2^2\big)\,\Big].
\end{eqnarray}

To the best of the authors' knowledge, there are no interior, a priori estimates of the gradient of the solution of (\ref{F:EllipticBVP}) (with $f=0$) as sharp as the bound provided by (\ref{Th:GandT}) for the solution itself. Therefore, based on more particular results such as \cite[Theorem 3.9]{G&T_book} and \cite[Problem 3.6]{G&T_book}, we make the reasonable assumption that for the PDE $(L+c)u=0$ in $\Omega_k$ with Dirichlet BCs on $\partial\Omega_k$
\begin{equation}
\label{F:Extended_GT}
\sup_{{\bf x}\in\Omega_k}||\nabla u({\bf x})||_2 \leq
K'_k+K''_k\sup_{{\bf x}\in\partial\Omega_k}|u({\bf x})|,
\end{equation}

where $K'_k$ and $K''_k$ are positive and may depend on anything but the value of the Dirichlet BC. Then, there exist positive constants $K',K''$ and $s$ such that
\begin{eqnarray}
\left.
\begin{array}{l}
\sum\limits_{k=1}^m ||{\bar w}_k({\bf x})||_2 \leq K' + \gamma_R K''\\
\sum\limits_{k=1}^m ||{\tilde w}_k({\bf x};{\bm\omegaup})||_2 \leq K' + 2\gamma_R K''\sup\limits_{1\leq j \leq s} {\cal N}(\omegaup_j)
\end{array}
\right\} \textrm{  if ${\bf x}\in\Omega$,}
\end{eqnarray}

The values $(K',K'')$ are the $(K'_k,K''_k)$ of the subdomain $k$ with the largest gradient estimate, and $s$ its number of nodes. The $\gamma_R$ shows up due to the interpolation along the interfaces; and $2\sup_j {\cal N}(\omegaup_j)= \sup_j |{\cal N}(\omegaup_j)|$. This leads to the bounds

\begin{equation}
V[\psi({\bm\omegaup})] \leq
V[{\bar\psi}] + 2\Lambda E[\int_0^{\tau} Y_t^2dt]\frac{1}{q^2}\big( K' + 2\gamma_R K''\sup\limits_{1\leq j \leq s} {\cal N}(\omegaup_j) \big)^2 =: V[{\bar\psi}] + {{\tilde v}({\bm\omegaup})},
\end{equation}
\begin{equation}
\label{F:V[barpsi]}
V[{\bar\psi}]\leq 2\Lambda E[\int_0^{\tau} Y_t^2dt](K'+\gamma_R K'')^2 =: {\bar v}.
\end{equation}

As the final preparatory step, let us calculate the noise-to-signal ratio (NSR)--defined as $NSR[\cdot]=\sqrt{V[\cdot]}/E[\cdot]$--of the variable ${\bar v}+{\tilde v}({\bm\omegaup})$ (\ref{F:NSR}):

\begin{equation}
\label{F:NSR}
\frac{\sqrt{V_{{\bm\omegaup}}[{\bar v}+{\tilde v}({\bm\omegaup})]}}{E_{{\bm\omegaup}}[{\bar v}+{\tilde v}({\bm\omegaup})]}=
\frac{\sqrt{V_{\bm\omegaup}\Big[\big(1+2\gamma_R\frac{\textrm{K'}}{\textrm{K''}}\sup\limits_{1\leq j \leq s}{\cal N}(\omegaup_j)\big)^2\Big]}}
{q^2\big(1+\gamma_R\frac{\textrm{K'}}{\textrm{K''}}\big)^2+E_{\bm\omegaup}\Big[ \big(1+2\gamma_R\frac{\textrm{K'}}{\textrm{K''}}\sup\limits_{1\leq j \leq s}{\cal N}(\omegaup_j)\big)^2 \Big]}.
\end{equation}

\begin{table}[h!]
\begin{footnotesize}
\[\begin{array}{|l|lllll||l|lllll|}
\multicolumn{1}{c}{} & \multicolumn{5}{c}{\textrm{nodes per subdomain ($s$)}} & \multicolumn{1}{c}{} & \multicolumn{5}{c}{\textrm{nodes per subdomain ($s$)}}\\
\hline
K'/K'' & 10 & 10^2 & 10^3 & 10^4 & 10^5 & K'/K'' & 10 & 10^2 & 10^3 & 10^4 & 10^5   \\
\hline
0  & .54 & .31 & .20 & .15 & .12      
      & 0  & .54  & .31  & .20   & .15  & .12        \\
10^{-2}    & .54  & .31  & .20   & .15  & .12  
		& 10^{-2}    & .54  & .31  & .21   & .15  & .12 \\
10^{-1}      & .51  & .29  & .20   & .15  & .12    
		& 10^{-1}      & .53  & .30  & .20   & .15  & .12 \\
1       & .30  & .21  & .15   & .12  & .10    
		& 1       & .40  & .25  & .18   & .13  & .11 \\
10    & .048  & .037  & .031   & .027  & .024 
		& 10 & .094  & .073  & .060   & .052  & .046 \\
10^2  & .0047  & .0035  & .0029   & .0025  & .0023   
		& 10^2 & .0094  & .0070  & .0059   & .0051  & .0046 \\
\hline   
\multicolumn{1}{c}{} & \multicolumn{5}{c}{\gamma_R=1} & \multicolumn{1}{c}{} & \multicolumn{5}{c}{\gamma_R=2}\\
\end{array}\]
\end{footnotesize}
\caption{Simulated NSR over $10^5$ realizations of formula (\ref{F:NSR}) for $q=2$. If $K'=0$ the NSR is independent of $\gamma_RK''$.}
\label{T:NSR}
\end{table}

We are finally in a position to precisely estate our claim and the supporting heuristic. The ultimate goal is to use the deterministic value $V[{\bar\psi}]$ instead of a random $V[\psi({\bm\omegaup})]$ yielded by one simulation (note that $V[{\bar\psi}]\leq E_{{\bm\omegaup}}\big[V[\psi({\bm\omegaup})]\big]$). Therefore, it is important that the ratio $\sqrt{V_{{\bm\omegaup}}\big[V[\psi({\bm\omegaup})]\big]}/E_{{\bm\omegaup}}\big[V[\psi({\bm\omegaup})]\big]$ be small so that $V[\psi({\bm\omegaup}=0)]:=V[{\bar\psi}]\approx V[\psi({\bm\omegaup})]$. This ratio is problem-dependent 
but, in order to provide a rough estimate, we substitute $V[\psi({\bm\omegaup})]$ by its upper bound ${\bar v}+{\tilde v}({\bm\omegaup})$. Then, we simulate the NSR of the bound (Table \ref{T:NSR}) for realistic values $\gamma_R=1,2$ and over a broad range of $K''/K'$ (which captures the effect of the $L$, $c$, the geometry $\Omega$, and the PDD partition $\{\Omega_k\}_{k=1}^m$), and the typical number of nodes per subdomain, $s$. Since the NSR of the proxy is negligible in most of the scenarios (and specially as $s$ grows), we argue that the same should hold for $V[\psi({\bm\omegaup})]$. Obviously, specific problems would allow for sharper estimates.\newline

\subsection*{H4: Loss of correlation due to discretization}\label{SS:H4}
Given $a_0$ and its corresponding timestep $h_0$ from (\ref{F:N_y_h}), the correlation between $\phi$ and $\xi(a;{\omegaup})$ is better than between their discretized counterparts $\phi_{h_0}$ and $\xi_{h_0}(a;{\omegaup})$. This leads to an overestimation of the predicted decrease of variance and has a significant effect, especially if $a$ and $a_0$ are comparable, or if $\rho[\phi,\xi]\gtrsim.99$. We invoke {\bf H1} to justify the following perturbative ansatz:
\begin{equation}
Cov[\phi_{h_0},\xi_{h_0}(a;{\bm\omegaup})]\approx Cov[\phi,\xi(a;{\bm\omegaup})]+B_{\Xi}\big( Cov[\phi_{h_0},\xi_{h_0}] \big),
\end{equation}

where $B_{\Xi}\big( Cov[\phi_{h_0},\xi_{h_0}] \big)$ is the discrete covariance bias.
Since $Cov[\phi_{h_0},\xi_{h_0}]=E[\phi_{h_0}\xi_{h_0}]-u_{h_0}E[\xi_{h_0}])]$ and $h_0$ is small enough, $B_{\Xi}\big( Cov[\phi_{h_0},\xi_{h_0}] \big)= {\cal O}(h_0^{\delta})$. That covariance bias does not seem accessible without having $\xi_{h_0}$, but we can use the fact that $Cov[\phi,\xi]=-V[\phi]$ (Lemma \ref{Th:Lema1}) to argue that
\begin{equation}
\Big|B_{\Xi}\big( Cov[\phi_{h_0},\xi_{h_0}] \big)\Big|\approx \Big|B_{\Xi}\big( V[\phi_{h_0}] \big)\Big|=:|\alpha|h^{\delta},
\end{equation} 
which can be extracted from a fit (see {\bf H5}). Then, assuming further that
\begin{equation}
\frac{1}{V[\phi_{h_0}]}\approx \frac{1}{V[\phi]} \qquad\textrm{   and   }\qquad
\frac{1}{V[\xi_{h_0}(a;{\bm\omegaup})]}\approx \frac{1}{V[\xi(a)]},
\end{equation} 

the squared correlation of the discretized variables is

\begin{equation}
\label{F:rho2_h}
\rho^2[\phi_{h_0},\xi_{h_0}(a;{\bm\omegaup})] \approx \rho^2[\phi,\xi(a;{\bm\omegaup})] - \frac{2\big|\alpha\rho[\phi,\xi(a;{\bm\omegaup})]\big|}{V[\phi]}h_0^{\delta},
\end{equation}

and since $|\beta| h_0^{\delta}=a_0/2$, the effective variance reduction in IterPDD($a_0,a$) is

\begin{equation}
\label{F:Final_formula}
\frac{V[\phi_{h_0}+\xi_{h_0}(a;{\bm\omegaup})]}{V[\phi_{h_0}]}\approx \frac{V[\phi+\xi(a;{\bm\omegaup})]}{V[\phi]} +
\big|\frac{\alpha\rho[\phi,\xi(a)]}{\beta V[\phi]} \big|a_0.
\end{equation}

For $a_0$ small enough, as $a\shortrightarrow a_0^+$, $|\rho[\phi,\xi(a,{\bm\omegaup})]|\shortrightarrow 1^-$, so (\ref{F:Final_formula}) is capped by 
\begin{eqnarray}
\label{F:VR_cap}
\lim_{a\shortrightarrow a_0^+} \frac{V[\phi_{h_0}+\xi_{h_0}(a;{\bm\omegaup})]}{V[\phi_{h_0}]} 
= \big| \frac{\alpha}{\beta V[\phi]} \big|a_0= \frac{2|B_{\Xi}(V[\phi])|}{V[\phi]}. 
\end{eqnarray}

This makes intuitively sense, because the variance of the score can hardly drop below its discretization error, even with an exact control variate.

\subsection*{H5: Fast estimation of constants}\label{SS:H5}

In order to apply the sensitivity formula (\ref{F:Sensitivity_formula}) and Algorithms \ref{A:IterPDD} and \ref{A:Scheduling}, a number of constants must be estimated. Here, we describe a fast way of accomplishing this. We stress the fact that any Monte Carlo simulation (whether or not related to PDD) also would require $\delta,\beta$ and $V[\phi]$ in order to enforce a set error tolerance. \newline


{\bf Global constants ($\delta$ and $\kappa$).} As already discussed, with smooth BVPs, integrators with at least approximately known $\delta$ should be available \cite{OurSurvey}. The constant $\kappa$ can easily be found by comparing the time taken by the computer to complete a number of visits to the implementations of (\ref{F:VR_Milstein_sys}) with and without ${\bf F}$.\newline 

{\bf Nodal constants related to first moments.} They are $E[\tau]$ (needed for $K_i$), $Cov[\phi,{\bar\psi}]$, and $\beta$. We consider a generic discretized random variable $\eta_h=\{\tau_h,{\bar\psi}_h\phi_h,\phi_h\}$, and assume that its first moment obeys the noisy model  
\begin{equation}
\label{F:Noisy_model_moment1}
E[\eta_h]\sim E[\eta_0] + {\cal N}(B^{(1)}h^{\delta},V[\eta_0]).
\end{equation}

Let $h_1>h_2>\ldots>h_{\hat M}$ be a 'cloud' of ${\hat M}$ equispaced timesteps. Set a number ${\hat N}$ and let ${\hat E}[\eta_h]$ be the mean of ${\hat N}$ realizations of $\eta_h$. After computing ${\hat M}$ MC independent simulations, the cloud of data
$\big(h_1,{\hat E}[\eta_{h_1}]\big),\ldots,\big(h_{\hat M},{\hat E}[\eta_{h_{\hat M}}]\big)$ is fitted to the noisy model (\ref{F:Noisy_model_moment1}) in order to extract $E[\eta_0]$ and $B^{(1)}$. (In order to do so, $V[\eta_0]$ in (\ref{F:Noisy_model_moment1}) can be replaced by the mean of the sample variances ${\hat V}[\eta_{h_1},\ldots,{\hat V}[\eta_{h_{\hat M}}]$.) 

Here, we provide a rougher recipe to carry out the fit with Matlab (see also \cite{LTMCR}). For this purpose, it is convenient to think of model (\ref{F:Noisy_model_moment1}) as a member of the generalized linear model (GLM) family
. Since $E[{\cal N}(B^{(1)}h^{\delta},V[\eta_0])]=B^{(1)}h^{\delta}$, the link is the identity and the fit is readily carried out by issuing the Matlab command
   
\begin{verbatim}
Coeff= glmfit((h.^delta),Eh,'normal','link','identity')
\end{verbatim}

where \verb|h| and \verb|Eh| above are Matlab arrays with respectively $(h_1,...,h_{\hat M}$ and $(\,{\hat E}[\eta_{h_1}],...,{\hat E}[\eta_{h_{\hat M}}]\,)$; and the components of the output, \verb|Coeff(1)| and \verb|Coeff(2)|, are the fitted values to $E[\eta_0]$ and $B^{(1)}$, respectively (check the Matlab documention for getting error bounds alongside).   




By applying this recipe, one gets approximations to the following quantities: if $\eta_h={\bar\psi}_h\phi_h$, to $Cov[\phi,{\bar\psi}]\approx {\hat Cov}[\phi,{\bar\psi}]= \verb|Coeff(1)|$ (and to its bias as a byproduct); and if $\eta_h=\phi_h$, to $\beta\approx {\hat\beta}=\verb|Coeff(2)|$ and to $E[\phi]\approx{\hat E}[\phi]=\verb|Coeff(1)|$ as a byproduct. The values ${\hat E}[\tau_{h_1}],...,{\hat E}[\tau_{h_{\hat M}}]$ are obtained along the means ${\hat E}[\phi_{h_1}],...,{\hat E}[\phi_{h_{\hat M}}]$ (see last line in Algorithm \ref{A:Gobet_Menozzi}), and are used to fit $E[\tau]$. \newline 

{\bf Nodal constants related to second moments.} Regarding the fitting of variances, it is well-known that variances of i.i.d. Gaussian distributions obey a scaled chi-squared PDF \cite{Kloeder&Platten99}. Accommodating the discretization bias, the appropriate noisy model is
\begin{equation}
\label{F:Noisy_model_moment2}
V[\eta_h]\sim B^{(2)}h^{\delta}
+\frac{V[\eta_0]}{{\hat N}-1}\chiup^2_{{\hat N}-1}=
B^{(2)}h^{\delta}
+\Gamma\big( \frac{{\hat N}-1}{2},\frac{2V[\eta_0]}{{\hat N}-1}\big)
\end{equation}
where $\chiup^2_{{\hat N}-1}$ is the Chi-squared distribution with ${\hat N}-1$ degrees of freedom, and $\Gamma(p_1,p_2)$ is a Gamma distribution with shape parameter $p_1=({\hat N}-1)/2$ and scale parameter $p_2=2V[\eta_0]/({\hat N}-1)$. Then, (\ref{F:Noisy_model_moment2}) can be identified with a GLM where the noise is Gamma and the link function is the identity, since 

\begin{equation}
E\big[V[\eta_h]\big]= B^{(2)}h^{\delta} + E\Big[\Gamma\big( \frac{{\hat N}-1}{2},\frac{2V[\eta_0]}{{\hat N}-1}\big)\Big]= B^{(2)}h^{\delta}+V[\eta_0].
\end{equation}




Let ${\hat V}[\eta_h]= \sum_{j=1}^{{\hat N}}(\eta_h^{(j)}-{\hat E}[\eta_h])^2\big/({\hat N}-1)$ be the sample variance. After filling the Matlab array \verb|Vh| with ${\hat V}[\phi_{h_1}],\ldots,{\hat V}[\phi_{h_{\hat M}}]$, issuing the command
\begin{verbatim}
Coeff= glmfit((h.^delta),Vh,'gamma','link','identity')
\end{verbatim}

yields \verb|Coeff(1)|$\approx V[\eta_0]/({\hat N}-1)$ and \verb|Coeff(2)|$\approx B^{(2)}$. Particularizing to $\eta_h=\phi_h$ allows to estimate ${\hat V}[\phi]\approx V[\phi]$ and ${\hat\alpha}\approx \alpha$; and to $\eta_h={\bar\psi}_h$ yields the fitted ${\hat V}[\phi]\approx V[{\bar\psi}]$ (and its bias).

\begin{algorithm}[h]
\caption{Fast fit of nodal constants}
\begin{algorithmic}[1]
\FOR{$k=1,\ldots,m$}
	\STATE{Solve ${\bar w}_k({\bf x})$ in (\ref{F:Error_propagation_function})}
\ENDFOR
\STATE{Store $\nabla {\bar w}({\bf x})=\oplus_{k=1}^m\nabla {\bar w}_k({\bf x})$}
\STATE{Estimate $\Pi$ and ${\tilde\Pi}$}
\FOR{$i=1,\ldots,n$}
   \STATE{{\bf set} ${\hat M},{\hat N},h_1,h_{{\hat M}}$ (may be point-dependent)}
	\FOR{$j=1,\ldots,{\hat M}$}
		\STATE{Run ${\hat N}$ realizations of Algorithm \ref{A:Gobet_Menozzi} at $h_j$, drawing ${\bar \psi}_h$ alongside}
		\STATE{Compute ${\hat E}[\phi_{h_j}],{\hat V}[\phi_{h_j}],{\hat E}[\tau_{h_j}],{\hat E}[{\bar\psi}_{h_j}],{\hat V}[{\bar\psi}_{h_j}\phi_{h_j}]$}
	\ENDFOR
	\STATE{Fit the nodal constants for this node}
\ENDFOR 
\STATE{{\bf output:} the fitted values $\{ {\hat E}_i[\phi], {\hat \beta}_i,{\hat V}_i[\phi],{\hat \alpha}_i,{\hat K}_i,{\hat \rho}_i[\phi,{\bar\psi}]\}_{i=1}^n$}
\end{algorithmic}
\label{A:Fitting}
\end{algorithm}

The full fitting procedure is sketched as Algorithm \ref{A:Fitting}. 
Note that the same sets of trajectories can be used for all the poinwise constants--it is only the functionals that change.

\subsection*{Sensitivity formula, scheduling-, and final algorithms}\label{SS:FinalAlgorithm}
Combining heuristics {\bf H1} through {\bf H5} we put forward the sensitivity formula (\ref{F:Sensitivity_formula}), based on (\ref{F:The_key}) and (\ref{F:Final_formula}) with ${\bar\psi}$ replacing $\psi({\bm\omegaup})$:

\begin{equation}
\label{F:Sensitivity_formula}
\rho^2[\phi_{h_0},\xi_{h_0}]\approx \rho^2[\phi,\xi]-\big|\frac{\alpha\rho[\phi,\xi]}{\beta}\big|a_0,
\textrm{  with }
\rho^2[\phi,\xi]\approx \frac{V[{\bar\psi}]a^2}{4V[\phi]}\big( 1-\rho^2[\phi,{\bar\psi}]\big).
\end{equation}

With it, a scheduling algorithm could be as Algorithm \ref{A:Scheduling}. There, the stopping criterion should take into account the quality of the available estimates and thus be problem-dependent. (We discuss this aspect further on Section \ref{S:Experiment}.)

\begin{algorithm}[h!]
\caption{Scheduling algorithm based on sensitivity formula (\ref{F:Sensitivity_formula}).}
\begin{algorithmic}[1]
\STATE{{\bf data}: $\delta,{\hat\kappa}$, and $\{{\hat K}_i,{\hat\beta}_i,{\hat\alpha}_i,{\hat V}_i[\phi],{\hat V}_i[{\bar\psi}],{\hat \rho}_i[\phi,{\bar\psi}]\}_{i=1}^n$}
\STATE{Set $j=0$, and $a_0$ is the nodal target error tolerance.}
\WHILE{stopping criterion not fulfilled}
\STATE{Solve the minimization problem
\begin{equation}
a_{j+1}= \arg\min\limits_{a>a_j} \sum\limits_{i=1}^n{\hat K}_i{\hat V}_i[\phi]\Bigg(
{\hat\kappa} \big( 1-{\hat r}_i^2(a) + \big|\frac{{\hat\alpha}_i{\hat r}_i(a)}{{\hat\beta}_i}\big|a_j \big) + \big(\frac{a_j}{a}\big)^{2+1/{\delta}}
\Bigg),
\end{equation}
with ${\hat r}_i^2(a)= 1- \frac{{\hat V}_i[{\bar\psi}]a^2}{4{\hat V}_i[\phi]}\big(1-{\hat\rho_i}^2[\phi,{\bar\psi}]\big)$.}
\STATE{Set $j=j+1$}
\ENDWHILE
\STATE{{\bf output:} $J=j$ and $a_J>a_{J-1}>\ldots>a_1>a_0$}
\end{algorithmic}
\label{A:Scheduling}
\end{algorithm}


Finally, we summarize the new version of PDD in Algorithm \ref{A:IterPDD}.

\begin{algorithm}[h!]
\caption{{\bf Iterative ('multigrid') probabilistic domain decomposition}}
\begin{algorithmic}[1]
\STATE{{\bf Data:} a global error tolerance $\epsilonup>0$ for a BVP like (\ref{F:EllipticBVP}) in $\Omega$, a confidence interval $P_q$}
\STATE{{\bf Choices:} PDD partition $(\{{\bf x}_i\}_{i=1}^n,\{\Omega_k\}_{k=1}^m)$, integrator $\Xi$ with known $\delta$, subdomain solver, interfacial interpolator $R$}
\STATE{Set the nodal target tolerance $a_0(\epsilonup)$ according to (\ref{F:a(eps)})}
\STATE{Solve $\{{\bar w}_k({\bf x})\}_{k=1}^m$ in (\ref{F:w_bar}) in parallel and construct ${\bar\psi}$ from Definition \ref{D:definition_of_psi} and (\ref{F:Bar_psi})}
\STATE{Estimate ${\hat\kappa}$ and the nodal constants with Algorithm \ref{A:Fitting}}
\STATE{Find $a_J>...>a_1$ with the scheduling algorithm (Algorithm \ref{A:Scheduling})}
\STATE{(Optional) Construct ${\tilde u}_J({\bf x})$ based on the 'fitted nodal values' ${\hat E}_1[\phi],\ldots,{\hat E}_n[\phi]$}
\FOR{$j=J...1$}
	\FOR{$i=1...n$}
		\STATE{For node ${\bf x}_i$ and $a_j$, determine $h$ and $N$ (\ref{F:N_y_h})}
		\STATE{Run $N$ independent realizations of Algorithm \ref{A:Gobet_Menozzi}}
		\STATE{Calculate the nodal value $u_i^{(j-1)}$}
	\ENDFOR
	\FOR{$k=1...m$}
\STATE{\begin{equation}
		\textrm{Solve the subdomain BVP }
		\left\{
		\begin{array}{ll}
		Lv_k + cv_k= f, & \textrm{ if } {\bf x}\in\Omega_k,\\
		v_k= g,			 & \textrm{ if } {\bf x}\in\partial\Omega_k\cap\partial\Omega,\\
		v_k= R[u_i^{(j-1)}|{\bf x}_i\in\Gamma_p^k], & \textrm{ if } {\bf x}\in\Gamma_p^k,\,1\leq p\leq m_k.
		\end{array}
		\right.\nonumber
		\end{equation}}
	\ENDFOR
	\STATE{Construct and store ${\tilde u}_{j-1}({\bf x})=\oplus_{k=1}^m v_k({\bf x})$ and $\nabla{\tilde u}_{j-1}({\bf x})=\oplus_{k=1}^m \nabla v_k({\bf x})$}
\ENDFOR
\end{algorithmic}
\label{A:IterPDD}
\end{algorithm}

\section{Numerical experiment} \label{S:Experiment}
\begin{figure}[h]
\centerline{\includegraphics[width=1.\columnwidth]{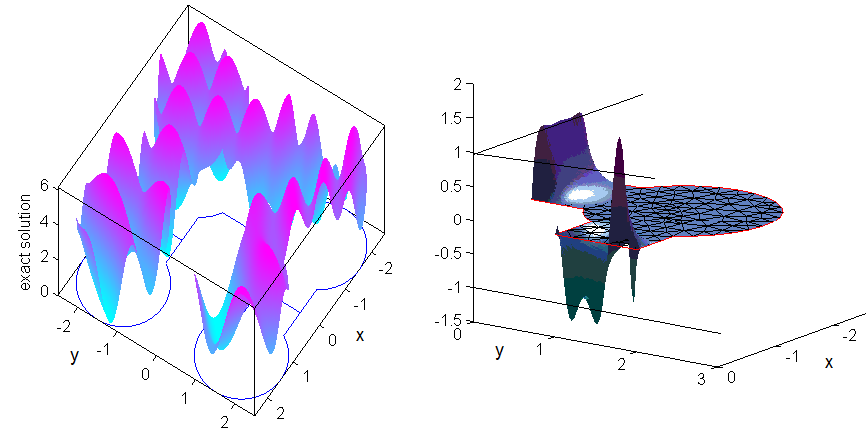}}
\caption{{\bf Left.} Exact solution of (\ref{F:PDE_ejemplo}). {\bf Right.} The bias-related function ${\bar w}_1(x,y)$ (\ref{F:Error_propagation_function}) used to generate ${\bar\psi}$ (\ref{F:definition_of_psi}) on $\Omega_1$ in Figure \ref{I:Figura1}. Notice the various signs of $\beta_i$ and the overshoots of the RBF interpolator. Here, $\gamma_R\approx 1.5$.}
\label{I:Figura3}
\end{figure}

We consider a BVP like (\ref{F:EllipticBVP}) with on the two-dimensional domain sketched in Figure \ref{I:Figura1} with $m=4$ subdomains, $3$ interfaces and $n=18$ nodes. The PDE is 
\begin{equation}
\label{F:PDE_ejemplo}
\nabla^2 u + \frac{\cos{(x+y)}}{1.1+\sin{(x+y)}}\big(\frac{\partial u}{\partial x}+\frac{\partial u}{\partial y}\big)-\frac{x^2+y^2}{1.1+\sin{(x+y)}}u+f(x,y)= 0,
\end{equation}
with $f(x,y)$ such that $u(x,y)=2\cos\big(2(y-2)x\big)+\sin\big(3(x-2)y\big)+3.1$ is the exact solution, as well as the Dirichlet BC--see Figure \ref{I:Figura3} (left). The coefficients of the stochastic representation of (\ref{F:PDE_ejemplo}) are: $\sigma=\sqrt{2}I_2$ ($I_2$ is the two-dimensional identity matrix), ${\bf b}= \cos{(x+y)}[1,1]\,/\big(1.1+\sin{(x+y)}\big)$, and $c=-(x^2+y^2)/\big(1.1+\sin{(x+y)}\big)<0$. The integrator is Algorithm \ref{A:Gobet_Menozzi}, for which we assume that $\delta=1$; $R$ is an RBF interpolator with multiquadrics \cite{Fasshauer_libro}, and the subdomain solver is FEM. The parameters of $R$ and FEM were so chosen that their errors are negligible compared with the nodal errors.\newline

{\bf Objectives.} With just four small subdomains, it obviously does not take a parallel computer to solve this toy problem--but it serves to test the idea, theory and approximations introduced in this paper in a controlled environment. For the sake of clarity, we stress that we will {\em not} compare the new IterPDD algorithm with results from deterministic domain decomposition, but with the previous version of PDD, PlainPDD. IterPDD inherits the scalability properties of PlainPDD, and comparisons of PlainPDD with deterministic domain decomposition methods can be found in \cite{Acebron2005,Acebron2009,Acebron2013} and references therein. Above all, we will focus on the speedup of IterPDD over PlainPDD, defined as
\begin{equation}
\label{F:S}
S(a_{j-1},a_j)= \frac{\textrm{cost of PlainPDD$(a_{j-1})$}}{\textrm{cost of IterPDD$(a_{j-1},a_j)$}}.
\end{equation}

The algorithms have been coded in fully vectorized Matlab, and the nodal values and subdomain BVPs solved on single processors sequentially, for ease of implementation. This suffices to study the computational costs and hence the speedup (\ref{F:S}). (The unit cost is a visit to Algorithm \ref{A:Gobet_Menozzi}). Particularly, we are interested in checking whether:
\begin{itemize}
\item the heuristics introduced in Section \ref{S:Approximations} are adequate;
\item the theoretically derived sensitivity formula (\ref{F:Sensitivity_formula}), correctly predicts the correlation of scores with pathwise control variates; and
\item the scheduling algorithm (Algorithm \ref{A:Scheduling}) gives a good approximation to the fastest sequence of IterPDD simulations yielding a final accuracy $a_0$.
\end{itemize}

\subsection{Preliminary numerical study of the speedup with perfectly balanced Monte Carlo simulations}\label{SS:Reference_sims}
Actual IterPDD simulations shown later will not be strictly balanced in the sense of Definition \ref{D:Balanced_MC_simulation} because they will rely of heuristic estimates {\bf H5} of the nodal constants according to Section \ref{SS:H5}. It is therefore informative to first carry out an exploration of the parameter space, using (nearly strictly) balanced Monte Carlo simulations. They will provide best-case reference values against which later, realistic IterPDD simulations in Section \ref{SS:Actual_sims} can be assessed.

\begin{table}[h]
\begin{footnotesize}
\[\begin{array}{|l|ll|llll|}
\multicolumn{1}{c}{\textrm{\bf Illustrative values}} & \multicolumn{6}{c}{\textrm{\bf{Accuracy $a_1$ of a previous rough PDD solution}}}\\
\hline
\multicolumn{1}{c}{} & \multicolumn{1}{l}{\textrm{exact}}&\multicolumn{1}{l}{\textrm{lookup}}&.02&.10&.26&\multicolumn{1}{l}{.62}\\ 
\hline 
\sum_{i=1}^nN_0^{(i)}\nu_0^{(i)}\textrm{ (with $\xi_0\big({\tilde u}_1(a_1)\big)\,$)}&  1.03\times 10^8	  & 1.01\times 10^8	&1.05\times 10^8 & 1.18\times 10^8 & 1.58\times 10^8 &5.82\times 10^8\\
V[\phi_{h_0}+\xi_{h_0}(a_1)] 	& 0.16	& 0.16 & 0.16 & 0.17 & 0.19 & 0.47\\
\hline       
\textrm{cost of PlainPDD$(a_1)$}&  & & 1.60\times 10^9 & 1.47\times 10^7 & 8.43\times 10^5 & 65989 \\
V[\phi_{h_1}]    & & & 4.92  & 5.10  & 5.45 & 4.08 \\
\big( \sum_{i=1}^n |u({\bf x}_i)-{\tilde u}_1({\bf x}_i)| \big)/n  & \multicolumn{2}{c|}{\textrm{use exact $u$ instead of ${\tilde u}(a_1)$}} & 0.012  & 0.06 & 0.11 & 0.77 \\
\sup_{\Omega}|u-{\tilde u}_1|  & & & 0.03 & 0.12 & 0.24 & 0.34 \\
\sup_{\Omega}||\nabla {\tilde u}_1||_2 & & & 0.87 & 0.94 & 2.14 & 6.61 \\
\hline
|\rho| & 0.9859 & 0.9859 & 0.9859 & 0.9851 & 0.9833 & 0.9595\\
S(a_0,a_1)   & 128.16 &72.61  & 7.37  & 58.14 & 46.39 & 22.67 \\
\hline                                                                                             

\hline
\end{array}\]
\end{footnotesize}
\caption{Acceleration by pathwise control variates, and other illustrative quantities, for reference (perfectly balanced) IterPDD$(a_0=.01,a_1)$ simulations.}
\label{T:Illustration}
\end{table}

As a preparatory step, we carefully compute $\kappa$
and the nodal constants (see Section \ref{SS:H5}). These values can be deemed exact and, in particular, the precise values of $\{\beta_i\}_{i=1}^{18}$ ensure that the pointwise MC simulations are balanced. Then, we pick $a_0=\{.0025,.005,.01,.02,.04\}$; $a_1=\{.02,.04,.06,.10,.14,\ldots,.62\}$; and run a set of $[{\tilde u}_1,\xi_0({\tilde u}_1)]=$IterPDD$(a_0,a_1)$ simulations that will serve as reference. For a target accuracy $a_0=.01$, Table \ref{T:Illustration} shows the dependence of some typical quantities with respect to the accuracy $a_1>a_0$ used to construct the rough ${\tilde u}(a_1)$ global solution. 
(Recall that $\xi_0$ is shorthand for $\xi_{h_0}$, etc.) 

On the column labeled 'exact', the control variate from the exact solution $\xi_0(u_{ex})$ has been used, so that all of the error is quadrature error. This represents the maximum benefit that could possibly be extracted from pathwise control variates. Otherwise, the approximation $\nabla {\tilde u}_1$ is stored on a lookup table (here, a grid), and interpolated from there to compute $\xi_0({\tilde u}_1)$--except on the column 'lookup', where the lookup table has been filled with the exact $\nabla u_{ex}$, in order to gauge purely the effect of interpolation. Since the cost of solving the subdomains and filling the lookup tables is negligible (i.e. $\Pi\approx 0 \approx {\tilde\Pi}$), the cost of IterPDD($a_0,a_1$) is essentially measured as:
\begin{equation}
\textrm{cost of IterPDD$(a_0,a_1)$}= \sum_{i=1}^{n} \big( N_1^{(i)}\nu_1^{(i)} + \kappa N_0^{(i)}\nu_0^{(i)}\big),
\end{equation}

where $N_0^{(i)}$ is the number of trajectories from node ${\bf x}_i$ at $h=h_0$, and so on.

Let us explain in detail the rightmost column of Table \ref{T:Illustration}. Without control variates, PlainPDD solves the BVP (\ref{F:PDE_ejemplo}) to within a nodal accuracy $a_0=.01$ at a cost $1.32\times 10^{10}$. If, instead, we get first a rougher PlainPDD solution (at $a_1=.62$), and use it to construct pathwise control variates, the total cost is $65989$ (for the rougher simulation, at $a_1$) plus $5.82\times 10^8$ (for the finer, at $a_0$); i.e. nearly $22.67$ times less. This is so because the average variance over the $n=18$ interfacial nodes has dropped, for the finer simulation at $h_0$, from about $4.90$ to $0.47$. Just how much the variance drops depends on the correlation $\rho[\phi_{h_0},\xi_{h_0}(a_1)]$; in this case the average correlation (in absolute value) is $|\rho|=0.9595$, which is quite good given that $a_1$ is $62$ times larger than $a_0$. The actual mean interfacial error, $0.77$, actually overshoots $a_1$; this can happen with a small probability, or if $h_1$ is not small enough. The global quality of the approximation ${\tilde u}_1={\tilde u}(a=a_1)$ is gauged by $\sup_{\Omega}|u-{\tilde u}_1|$. The connection between $\sup_{\Omega}|u-{\tilde u}_1|$ and $\sup_{\Omega}||\nabla{\tilde u}_1||_2$ with the mean interfacial error were discussed in Section \ref{S:GlobalError} and {\bf H5}, and the values in this concrete example are given here for the sake of illustration.



\begin{figure}[h!]
\centerline{\includegraphics[width=1.\columnwidth]{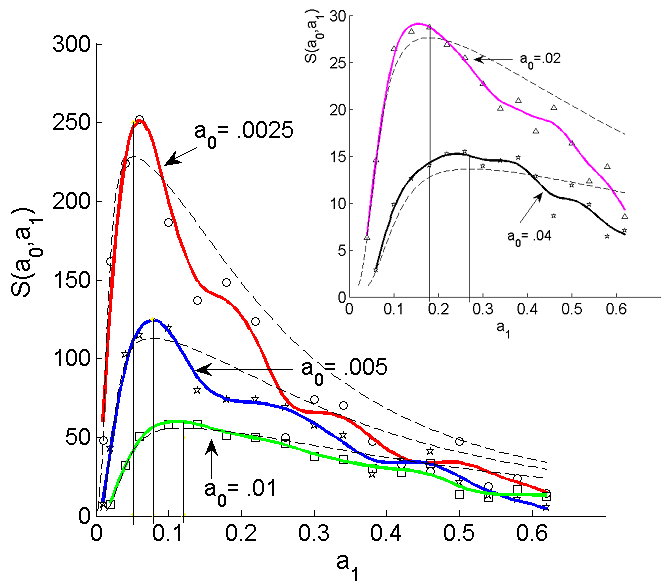}}
\caption{Speedup of IterPDD$(a_0,a_1)$ over PlainPDD$(a_0)$. The symbols denote the reference (strictly balanced) PDD simulations. In order to show the trend, they have been approximated with smoothing splines (solid color curves). {\bf Dashed curves}: predicted speedup, based on the sensitivity formula (\ref{F:Sensitivity_formula}) and fast estimates of constants (Algorithm \ref{A:Fitting}). Although the latter do not lead to balanced MC simulations as with the reference simulations, their maxima calculated with Algorithm \ref{A:Scheduling} (vertical lines) are coincident--see (*) on Table \ref{T:Schedula}.}
\label{I:Figura4}
\end{figure}

On Figure \ref{I:Figura4}, the values of  $S(a_0,a_1)$ for the full set of reference simulations are depicted with symbols. Recall that those speedups are realizations of some underlying distribution, labeled with ${\bm\omegaup}$ (see Section \ref{SS:H3}). Nonetheless, as a guide for the eye we have joined those symbols corresponding to the same target accuracy $a_0$ with smoothing splines. As $a_1$ grows, the effect of randomness on $S(a_0,a_1)$ is more significant.

The symbols and solid lines on Figure \ref{I:Figura4} summarize our preliminary numerical investigation. In the most favourable test (when $a_0=.0025$ and $a_1\approx 0.05$, i.e. the peak of the red curve), a speedup of $250$ times was achieved. Of course, this maximum is obtained {\em a posteriori} and relying of perfect information (the exact nodal constants). The object of this paper is to {\em predict} the maxima on Figure \ref{I:Figura4}. In what follows, we will show that i) the optimal $a_1$ for each target accuracy $a_0$ can be quite accurately predicted by the sensitivity formula (\ref{F:Sensitivity_formula}); ii) the speedup itself can also be well predicted; and iii) it is possible to improve further the speedup by iterating the method, as dictated by the scheduling algorithm (Algorithm \ref{A:Scheduling}). Those predictions rely on the fast estimates of nodal constants by heuristic {\bf H5} (Section \ref{SS:H5}), as well as on the auxiliar global variable ${\bar\psi}$. Critically, all the needed ingredients can be obtained at a relatively negligible cost and in parallel.

\subsection{Heuristics {\bf H1-H5} and actual IterPDD simulations}\label{SS:Actual_sims}

Now, we proceed to test Algorithm \ref{A:IterPDD} proper. First, we perform fast estimates of the nodal constants with Algorithm \ref{A:Fitting}, taking (without trying to optimize in any sense), ${\hat M}=100$ equispaced timesteps $h$ in $[.001,.01]$ and ${\hat N}=1000$. (As intended, the cost of this is comparatively very small.) A good estimate ${\hat\kappa}=1.8$ is straightforward to produce. Moreover, we solve (\ref{F:w_bar}) on each subdomain to construct ${\bar\psi}$ by Definition \ref{D:definition_of_psi} and by (\ref{F:Bar_psi})--see Figure \ref{I:Figura3} (right). Based on the resulting fitted constants, we run the scheduling Algorithm \ref{A:Scheduling} in order to get the 'optimal sequences' $a_J>a_{J-1}>\ldots>a_1>a_0$ for the same set $a_0=\{.0025,.005,.01,.02,.04\}$ as with the reference simulations before. The minimization (line $4$ in Algorithm \ref{A:Scheduling}) is carried out with Matlab's \verb|fmincon|, and stopped as soon as the predicted $S(a_j,a_{j+1})<1.5$. The calculated optimal sequence for each $a_0$ are given on the left section of Table \ref{T:Schedula}. Several illustrative quantities, both those predicted by Algorithm \ref{A:Scheduling}, and those actually attained after the simulation (using the same set of $\{{\hat \beta}_i\}_{i=1}^{18}$) are listed on the central and right sections of Table \ref{T:Schedula}, respectively. 

Let us explain the top vertical data block in Table \ref{T:Schedula}. We wish to get a PDD global solution with nodal accuracy $a_0=.04$. Algorithm \ref{A:Scheduling} calculates that the fastest way of getting it is by running three PDD iterations altogether, the last two of which use control variates. First, PlainPDD($a_2=.92$), next IterPDD($a_1=.27,a_2$), and finally IterPDD($a_0=.04,a_1$). It predicts that the total cost of this sequence is $13.93$ times less than directly computing PlainPDD($a_0$)--which directly translates into being $13.93$ times faster, given the embarrassingly parallel quality of PDD. This is called 'cumulative speedup', i.e. the cost of IterPDD$(a_0,a_1,\ldots,a_J)$ over that of PlainPDD$(a_0)$. It also predicts that there is no gain in getting one further previous iteration (at some $a_3>a_2$), because already $S(a_2,a_1)=1.95$ is close to the set threshold $1.5$. Compare also some of the values predicted by Algorithm \ref{A:Scheduling} with those actually observed {\em a posteriori}, on the two rightmost columns of Table \ref{T:Schedula}. Let us now move to the inset in Figure \ref{I:Figura4}. In the reference simulations reported in Section \ref{SS:Reference_sims}, the experimental curve for $S(a_0=.04,a_1)$ is represented by the smoothing spline in solid black, whose peak takes place at about $a_1=0.25$. Algorithm \ref{A:Scheduling} works with a model of that spline depicted by the dashed black curve hovering near it. The maximum of that approximation takes place at $a_1=.27$ (vertical line), which is very close, both in terms of the position and magnitude of the maximum. Moreover, Algorithm \ref{A:Scheduling} decides that it can still improve on that by running a previous PlainPDD($a_2=.92$) simulation, and exploiting the resulting control variates.

This is repeated for $a_0=\{.0025,.005,.01,.02,.04\}$. First, predictions made by Algorithm \ref{A:Scheduling} are shown on Table \ref{T:Schedula}--and compared with actually observed quantities. Second, the speedup curves predicted by Algorithm \ref{A:Scheduling} are also plotted with dashed lines on Figure \ref{I:Figura4}, and their maxima highlighted with vertical lines. A few comments are in order. 

\begin{table}[h!]
\begin{footnotesize}
\[\begin{array}{|l|lll|ll|}
\multicolumn{1}{c}{} & \multicolumn{3}{c}{\textrm{{\bf predicted with Algorithm \ref{A:Scheduling}}}} & \multicolumn{2}{c}{\textrm{{\bf observed with Algorithm \ref{A:IterPDD}}}}\\
\hline
\multicolumn{1}{l}{\textrm{{\bf 'optimal sequence'}}} & \multicolumn{3}{l|}{\textrm{{\bf cost of}}} & \multicolumn{1}{l}{\textrm{{\bf cost of}}} &\multicolumn{1}{c}{}\\
\multicolumn{1}{l}{\bm{ a_J\shortrightarrow\ldots\shortrightarrow a_1\shortrightarrow a_0}} & \multicolumn{1}{l}{\textrm{{\bf IterPDD$(a_{j-1},a_j)$}}} & \bm{|\rho|} & \bm{ S(a_{j-1},a_j)} & \textrm{{\bf IterPDD$(a_{j-1},a_j)$}} & \multicolumn{1}{l}{\bm{|\rho|}}  \\
\hline
a_2=.92\shortrightarrow a_1 & 3.38\times 10^5 & 0.677 & 1.95 & 3.64\times 10^5 & 0.830  \\
a_1=.27(*)\shortrightarrow a_0 & 1.54\times 10^7 & 0.960 & 13.63 & 1.49 \times 10^7& 0.970   \\
\hline
\multicolumn{1}{l}{{\bf a_0=.04}}& \multicolumn{3}{c|}{\textrm{predicted cumulative speedup$=13.93$}} &
       \multicolumn{2}{c}{} \\ 
\multicolumn{2}{c}{}&&&\multicolumn{1}{c}{}\\
\hline
a_2=.69\shortrightarrow a_1 & 7.84\times 10^5 & 0.776 & 2.97 & 9.09\times 10^5 & 0.864  \\
a_1=.18(*)\shortrightarrow a_0 & 6.07\times 10^7 & 0.981 & 27.62 & 5.41\times 10^7 & 0.985  \\
\hline
\multicolumn{1}{l}{{\bf a_0=.02}}& \multicolumn{3}{c|}{\textrm{predicted cumulative speedup$=28.34$}} & \multicolumn{2}{c}{} \\ 
\multicolumn{2}{c}{}&&&\multicolumn{1}{c}{}\\
\hline
a_2=.53\shortrightarrow a_1 & 1.80\times 10^6 & 0.862 & 4.54 & 2.06\times 10^6 & 0.915  \\
a_1=.12(*)\shortrightarrow a_0 & 2.40\times 10^8 & 0.991 & 55.90 & 2.19\times 10^8 & 0.992  \\
\hline
\multicolumn{1}{l}{{\bf a_0=.01}}& \multicolumn{3}{c|}{\textrm{predicted cumulative speedup$=57.42$}} & \multicolumn{2}{c}{}  \\ 
\multicolumn{2}{c}{}&&&\multicolumn{1}{c}{}\\
\hline
a_2=.41\shortrightarrow a_1 & 4.12\times 10^6 & 0.918 & 6.93 & 4.07\times 10^6 & 0.942 \\
a_1=.078(*)\shortrightarrow a_0 & 9.49\times 10^8 & 0.995 & 113.01 & 8.80 \times 10^8& 0.996 \\
\hline
\multicolumn{1}{l}{{\bf a_0=.005}}& \multicolumn{3}{c|}{\textrm{predicted cumulative speedup$=116.00$}} &\multicolumn{2}{c}{} \\
\multicolumn{2}{c}{}&&&\multicolumn{1}{c}{}\\
\hline
a_3=1.03\shortrightarrow a_2 & 2.49\times 10^5& 0.657 & 1.68 & 3.18\times 10^5 & 0.776   \\
a_2=.32\shortrightarrow a_1 & 9.40\times 10^6 & 0.948 & 10.59 & 8.78\times 10^6& 0.962   \\
a_1=.051(*)\shortrightarrow a_0 & 3.76\times 10^9 & 0.9977 & 228.22 & 3.94\times 10^9 & 0.9980 \\
\hline
\multicolumn{1}{l}{{\bf a_0=.0025}}& \multicolumn{3}{c|}{\textrm{predicted cumulative speedup$=233.84$}} &\multicolumn{2}{c}{} \\
\end{array}\]
\end{footnotesize}
\caption{Comparison of results predicted by the 'optimal sequence' obtained with the scheduling algorithm and those actually obtained after running Algorithm \ref{A:IterPDD} with the constants fitted with Algorithm \ref{A:Fitting}. See Figure \ref{I:Figura4} for the approximate maxima marked with (*), and text for further details.}
\label{T:Schedula}
\end{table}

Even though Algorithm \ref{A:Scheduling} relies on fast estimates of the constants, non-balanced simulations, and neglects the effects of randomness, the predicted speedup curves resemble quite well those obtained from the reference simulations. In particular, the position of the maxima of both sets are nearly coincident. Moreover, for $a_0,a_1$ small enough, Algorithm \ref{A:Scheduling} provides estimates of the costs, speedup and mean correlation (on Figure \ref{I:Figura4} and Table \ref{T:Schedula}) which are consistently conservative, as predicted by the theory in Section \ref{S:Approximations}. On the other side, when the tolerances $a_{j-1},a_j$ are not so small, heuristics {\bf H1-H4} break down, and the fact that $V[{\bar\psi}]<E_{{\bm\omegaup}}\big[V[\psi({\bm\omegaup})]\big]$ begins to tell, so the estimates may not be conservative, but still they are acceptable. At least with example (\ref{F:PDE_ejemplo}), most of the cumulative speedup is attained on the last IterPDD simulation. 


Finally, let us mention that running the scheduling algorithm with nodal constants based on $\psi({\bm\omega})$ instead of ${\bar\psi}$ (as we did) leads to very similar results to those reported. In fact, in this BVP, it seems that the effect of discretization on $\rho[\phi_h,\xi_h]$ (\ref{F:VR_cap}) outweighs the effect of randomness in the range of $a_0$ studied. This is why the speedups in Table \ref{T:Schedula} grow linearly with $a_0$. Importantly, this means that one order of magnitude has been knocked down from the Monte Carlo cost estimate ${\cal O}(1/a_0^{2+1/\delta})$ given in the Introduction,  to  ${\cal O}(a_0^{-2})$. 

\section{Conclusions} \label{S:Conclusions}
In this paper we have laid out the theoretical foundations of a much improved version of PDD, which we have called IterPDD.
The theoretical formulas have been derived for linear, smooth, second order elliptic BVPs with Dirichlet BCs. For this case, all the required ingredients of IterPDD are currently available: the probabilistic representation (Dynkin's formula), efficient SDE numerics (the Gobet-Menozzi integrator), and the Green's function representation of the solution. As long as those three items are in place, IterPDD can be extended to other BVPs (at least, linear ones)--although the specific formulas need to be adjusted correspondingly.  

With this goal, the PDD programme is currently being further developed in three main directions: i) blending it with Giles' Multilevel method \cite{Giles_y_yo}; ii) extending it to parabolic BVPs and mixed BCs; and iii) adjusting for processors with insufficient memory.

On the other hand, hyperbolic problems are typically difficult to handle with stochastic approaches, although representations for some of them do exist \cite{Acebron2016}. Elliptic and parabolic problems with discontinuous coefficients or boundary singularities can be handled in some cases \cite{Chip}, but often need specific stochastic representations and/or tailored numerical methods, which may not yet be satisfactory. Finally, while mildly nonlinear BVPs could be accommodated into the current IterPDD framework by linearization, of far greater interest are probabilistic representations of nonlinear equations, such as in \cite{Acebron2009,Acebron2013}. 

The purpose of this mainly theoretical paper is to introduce the strategy of pathwise control variates into the framework of probabilistic domain decomposition, as well as heuristics and approximations which make the idea useful in practice. The numerical results on the
example used for illustration are very encouraging, both in terms of confirming the heuristics and of performing hundreds of times faster than the previous version of PDD. It is clear, however, that substantially larger and more challenging problems must be tackled in order to assess the real potential of probabilistic domain decomposition. 

\section{Acknowledgements}
This work was supported by Portuguese national funds through FCT under grants UID/CEC/50021/2013 and PTDC/EIA-CCO/098910/2008. FB also acknowledges FCT funding under grant SFRH/BPD/79986/2011.



\begin{thebibliography}{10}
\small
\bibitem{Acebron2005}
 J.A. Acebr\'on, M.P. Busico, P. Lanucara and R. Spigler, {\em Domain
decomposition solution of elliptic problems via probabilistic
methods}. SIAM J. Sci. Comput. {\bf 27}, 440-457 (2005). 

\bibitem{Acebron2009}
J.A. Acebr\'on, A. Rodr\'iguez-Rozas and R. Spigler, {\em Domain decomposition solution of nonlinear two-dimensional parabolic problems by random trees}. J. Comput. Phys. {\bf 228}, 5574-5591 (2009).

\bibitem{Acebron2013}
J.A. Acebr\'on and A. Rodr\'iguez-Rozas, {\em Highly efficient numerical algorithm based on random trees for accelerating parallel Vlasov-Poisson simulations}. J. Comput. Phys. {\bf 250}, 224-245 (2013).

\bibitem{Acebron2016}
J.A. Acebr\'on and M. Ribeiro, {\em A Monte Carlo method for solving the one-dimensional telegraph equations with boundary conditions}. J. Comput. Phys. {\bf 305}, 29-43 (2016).

\bibitem{Chip}
F. Bernal, J.A. Acebr\'on and I. Anjam, {\em A Stochastic Algorithm Based on Fast Marching for Automatic Capacitance Extraction in Non-Manhattan Geometries}.
SIAM Journal on Imaging Sciences {\bf 7}(4), 2657-2674 (2014).

\bibitem{OurSurvey}
F. Bernal and J.A. Acebr\'on, {\em A Comparison of Higher-Order Weak Numerical Schemes for Stopped Stochastic Differential Equations}. Submitted (2015).

\bibitem{Bers1964}
L. Bers, F. John and M. Schechter, {\em Partial differential equations}. Interscience (1964)

\bibitem{Fasshauer_libro}
G.E. Fasshauer, {\em Meshfree Approximation Methods with Matlab}. Interdisciplinary Mathematical Sciences--Vol. 6 World Scientific Publishers, Singapore (2007).

\bibitem{Freidlin_Book}
M. Freidlin, {\em Functional integration and partial differential equations}. Annals of Mathematics Studies, vol. 109, Princeton University Press, Princeton (1985).

\bibitem{G&T_book}
D. Gilbarg and N.S. Trudinger, {\em Elliptic Partial Differential Equations of Second Order}. Revised $3^{rd}$ edition, Springer-Verlag Heidelberg Berlin (2001).

\bibitem{Giles_y_yo}
M.B. Giles and F. Bernal, {\em Multilevel simulations of expected exit times and other functionals of stopped diffusions}. In preparation (2015).

\bibitem{Glasserman2004}
P. Glasserman, {\em Monte Carlo Methods in Financial Engineering}. Springer-Verlag, New York (2004).


\bibitem{Gobet&Menozzi_2010}
E. Gobet and S. Menozzi, {\em Stopped diffusion processes: overshoots and boundary correction}. Stochastic Processes and their Applications, {\bf 120}, 130-162, (2010).

\bibitem{Haji-Ali}
A.L. Haji-Ali, F. Nobile, E. von Schwerin and R. Tempone, {\em Optimization
of mesh hierarchies in multilevel Monte Carlo samplers}, ArXiv preprint:
1403.2480 (2014)


\bibitem{Kloeder&Platten99}
P.E. Kloeden and E. Platten, {\em Numerical Solution of Stochastic Differential Equations}. Springer, Applications of Mathematics 23 (1999).

\bibitem{LTMCR}
S. Mancini, F. Bernal and J.A. Acebr\'on, {\em An efficient algorithm for accelerating Monte Carlo approximations of the solution to boundary value
problems}. Accepted in the Journal of Scientific Computing (2015).

\bibitem{Milstein&Tretyakov_book}
G. N. Milstein and M. V. Tretyakov, {\em Stochastic Numerics for Mathematical Physics}, Springer-Verlag, Berlin, 2004.

\bibitem{Smith96}
B.F. Smith, P.E. Bjorstad and W.D. Grop, {\em Domain Decomposition: Parallel Multilevel Methods for Elliptic Partial Differential Equations}. Cambridge University Press, Cambridge (1996).

\end{thebibliography}
\end{document}